\documentclass{amsart}
\usepackage{amsmath}
\usepackage{amsthm}
\usepackage{fullpage}   
\usepackage{pslatex}    

\newtheorem{maintheorem}{Theorem}
 
\newtheorem{theorem}{Theorem}[section]
\newtheorem{lemma}[theorem]{Lemma}
\newtheorem{proposition}[theorem]{Proposition}
\newtheorem{corollary}[theorem]{Corollary}
\theoremstyle{definition}
\newtheorem{definition}[theorem]{Definition}
\newtheorem{question}[theorem]{Question}
\newtheorem{remark}[theorem]{Remark}
\newtheorem{solution}[theorem]{Solution}
\numberwithin{equation}{theorem}
\newcommand{\SI}[1]{\textbf{SI(#1)}}
\newcommand{\MI}[1]{\textbf{MI(#1)}}
\newcommand{\boplus}{\textstyle \bigoplus}
\newcommand{\botimes}{\textstyle \bigotimes}
\newcommand{\type}[1]{{\ensuremath(#1)}}
\DeclareMathOperator{\Proj}{Proj}
\DeclareMathOperator{\Mod}{Mod}
\DeclareMathOperator{\adeg}{adeg}
\DeclareMathOperator{\Ext}{Ext}
\DeclareMathOperator{\Tor}{Tor}
\DeclareMathOperator{\End}{End}
\DeclareMathOperator{\cd}{cd}
\DeclareMathOperator{\Hom}{Hom}
\DeclareMathOperator{\HB}{H}
\DeclareMathOperator{\QGr}{QGr}
\DeclareMathOperator{\GKdim}{GKdim}

\DeclareMathOperator{\im}{im}
\newcommand{\fm}{\mathfrak{m}}
\newcommand{\Z}{\mathbb{Z}}
\begin{document}

\title[Regular algebras of dimension 4]
{Regular algebras of dimension 4 and \\
their $A_\infty$-Ext-algebras}

\author{D.-M. Lu, J. H. Palmieri, Q.-S. Wu and J. J. Zhang}

\address{(Lu) Department of Mathematics, Zhejiang University,
Hangzhou 310027, China}

\email{dmlu@zju.edu.cn}

\address{(Palmieri) Department of Mathematics, Box 354350, University of 
Washington, Seattle, WA 98195, USA}

\email{palmieri@math.washington.edu}

\address{(Wu) Institute of Mathematics, Fudan University, 
Shanghai, 200433, China}

\email{qswu@fudan.edu.cn}

\address{(Zhang) Department of Mathematics, Box 354350, University of 
Washington, Seattle, WA 98195, USA}

\email{zhang@math.washington.edu}

\begin{abstract}
We construct four families of Artin-Schelter regular algebras of global
dimension four. Under some generic conditions, this is a complete list
of Artin-Schelter regular algebras of global dimension four that are 
generated by two elements of degree 1. These algebras are also strongly 
noetherian, Auslander regular and Cohen-Macaulay. One of the main tools 
is Keller's higher-multiplication theorem on $A_\infty$-Ext-algebras.
\end{abstract}




\maketitle



\section*{Introduction}
\label{yysect0}

Recent developments in noncommutative projective algebraic geometry 
and its applications to fields such as mathematical physics 
demand to have more examples of quantum spaces. For example, 
what are the quantum K3 surfaces and the quantum Calabi-Yau 3-folds?
One naive idea is to construct these quantum spaces as subschemes 
of some higher dimensional quantum ${\mathbb P}^n$s --- noncommutative 
analogues of projective $n$-spaces. The quantum
${\mathbb P}^2$s have been classified by Artin, Schelter, Tate and 
Van den Bergh \cite{ASc, ATV1, ATV2} and understood very well by
many researchers. One of the central questions in noncommutative 
projective geometry is 

\medskip

\centerline{\it the classification of quantum ${\mathbb P}^3$s.}

\medskip

The complete classification of quantum ${\mathbb P}^3$s is an 
extremely difficult project and is probably an unreachable 
goal in the near future. An algebraic approach of constructing quantum 
${\mathbb P}^n$s is to form the noncommutative scheme $\Proj A$ where 
$A$ is a noetherian Artin-Schelter regular connected graded algebra 
of global dimension $n+1$. Therefore the algebraic version of the 
above mentioned question is

\medskip

{\it the classification of noetherian, Artin-Schelter regular, 
connected graded algebras of global dimension 4}. 

\medskip

Researchers have been studying many special classes of noetherian 
Artin-Schelter regular algebras of global dimension 4. For simplicity 
we only consider
graded algebras that are generated in degree 1. The most famous one is the 
Sklyanin algebra of dimension 4, introduced by Sklyanin \cite{Sk1, Sk2}. 
Homological and ring-theoretic properties of the Sklyanin algebra were 
understood by Smith and Stafford \cite{SS}.  Levasseur and Smith \cite{LS} 
studied the representations of the Sklyanin algebra. Some deformations 
of Sklyanin algebras were studied by Stafford \cite{Staf}. The homogenized 
$U({\mathfrak {sl}}_2)$ was studied by Le Bruyn and Smith \cite{LBS}.
Normal extensions of Artin-Schelter regular algebras of dimension 3 
were studied by Le Bruyn, Smith, and Van den  Bergh \cite{LBSV}. The quantum 
$2\times 2$-matrix algebra was studied by Vancliff \cite{Va1, Va2}. 
Some classes of Artin-Schelter regular algebras with finitely many 
points and Artin-Schelter regular algebras containing a 
commutative quadric were studied by Shelton, Van Rompay, Vancliff, 
Willaert, etc. \cite{SV1, SV2, VV1, VV2, VVW}. This is a partial 
list, and the different families may overlap. However, up to now, 
we do not have a clear picture of the complete classification. 

Note that all of the regular algebras listed in the previous paragraph are 
Koszul and generated by four elements of degree 1. The Koszul duals of 
some of these algebras have been studied also. From the list it seems 
that Koszul regular algebras are more popular than the non-Koszul 
ones. There are two explanations. One is that the non-Koszul regular 
algebras can be more difficult to study since the relations of such 
algebras are not quadratic. The other is that a non-Koszul algebra 
$A$ is not a deformation of the commutative polynomial ring 
$k[x_0,x_1,x_2,x_3]$ where $k$ is a base field, hence
$\operatorname{Proj}\; A$ is not a ``classical deformation'' of the 
commutative ${\mathbb P}^3$. Nevertheless it is important to search 
for both the Koszul and non-Koszul regular algebras of dimension four
or higher. 

Throughout, let $k$ be a commutative base field; in Section
\ref{yysect5}, we will assume that $k$ is algebraically closed.  Our
main result is the following.

\begin{maintheorem}
\label{yythm0.1} The following algebras are Artin-Schelter 
regular of global dimension four.
\begin{enumerate}
\item 
$A(p):=k\langle x,y\rangle/(xy^2-p^2y^2x,x^3y+px^2yx+p^2xyx^2+p^3 yx^3)$,
where $0\neq p\in k$.
\item
$B(p):=k\langle x,y\rangle/(xy^2+ip^2y^2x,x^3y+px^2yx+p^2xyx^2+p^3 yx^3)$,
where $0\neq p\in k$ and $i^2=-1$.
\item
$C(p):=k\langle x,y\rangle/(xy^2+ p yxy+p^2y^2x,x^3y+j p^3 yx^3)$,
where $0\neq p \in k$ and $j^2-j+1=0$.
\item
$D(v,p):=k\langle x,y\rangle/(xy^2+vyxy+p^2y^2x,x^3y+(v+p)x^2yx+
(p^2+pv)xyx^2+p^3 yx^3)$,
where $v,p\in k$ and $p\neq 0$.
\end{enumerate}
If $k$ is algebraically closed, then this list 
(after deleting some special algebras in each family) is, up to 
isomorphism, a complete list of $(m_2,m_3)$-generic Artin-Schelter 
regular algebras of global dimension four that are generated by two 
elements.
\end{maintheorem}

The $(m_2,m_3)$-generic condition is a condition on the 
multiplications $m_2$ and $m_3$ of the 
$\Ext$-algebras.  The details will be 
explained later---see (GM2) and (GM3) in Section \ref{yysect5}. The 
enveloping algebra of $U(\mathfrak g)$ given in \cite[1.21(iii)]{ASc} 
is a special case of (d) when $p=1$ and $v=-2$. None of the  
algebras in Theorem \ref{yythm0.1} is Koszul. 

Note that every algebra in Theorem \ref{yythm0.1} is 
${\mathbb Z}^2$-graded with $\deg x=(1,0)$ and $\deg y=(0,1)$.
The proof of Theorem \ref{yythm0.1} also shows:

\begin{maintheorem}
\label{yythm0.2} The list given in Theorem \ref{yythm0.1} is, up
to isomorphism, a complete list of noetherian Artin-Schelter regular 
${\mathbb Z}^2$-graded algebras of global dimension four generated 
by two elements with degrees $(1,0)$ and $(0,1)$.
\end{maintheorem}

Using the Artin-Schelter regular algebras given above, we can form 
new examples of quantum spaces of dimension three. We also show that 
algebras in Theorem \ref{yythm0.1} have good homological properties, 
which are useful for understanding the structure of their projective
schemes. 

\begin{maintheorem}
\label{yythm0.3} All of the algebras in Theorem \ref{yythm0.1} are strongly 
noetherian, Auslander regular and Cohen-Macaulay. 
\end{maintheorem}

A crucial step in proving Theorem \ref{yythm0.1} is to study the 
$A_\infty$-structure on Ext-algebras. The use of $A_\infty$-algebras 
in graded ring theory is a new approach since the concept 
of an $A_\infty$-algebra 
was introduced to noncommutative algebra very recently. This method has 
advantages for non-Koszul regular algebras because one can get 
information from the non-trivial higher multiplications on the 
Ext-algebras. It is plausible that all non-Koszul Artin-Schelter 
regular algebra of global dimension 4 will be classified by this 
method. A complete classification of 4-dimensional Artin-Schelter 
regular algebras generated by two elements of degree 1 is undergoing 
in \cite{LWZ}. This method could also be useful for constructing some 
new families of Artin-Schelter regular algebras of higher global 
dimension. 

Here is an outline of the paper.  In Section~\ref{yysect1}, we recall
the definition of Artin-Schelter regular algebras; we focus on those
which have global dimension 4 and are generated in degree 1, computing
their Hilbert series and the size of their Ext algebras.  There are
three types of such Artin-Schelter regular algebras, called \type{12221}, \type{13431},
and \type{14641}.  Those of type \type{14641} are Koszul, while the $(m_{2},
m_{3})$-generic algebras of Theorem~\ref{yythm0.1} are all of type
\type{12221}.

The goals of Section~\ref{yysect2} are to define $A_\infty$-algebras
and to state a theorem of Keller's relating the $A_{\infty}$-structure
on $\Ext_{A}^{*}(k_{A},k_{A})$ to the relations in $A$.  Some grading
issues are also discussed there.  In Section~\ref{yysect3}, we apply
the computations of the sizes of Ext algebras in Section~\ref{yysect1}
to compute the possible $A_\infty$-structures on Ext algebras of
regular algebras of dimension 4.  We specialize to algebras of type
\type{12221} starting in Section~\ref{yysect4}, and we introduce the
generic conditions in Section~\ref{yysect5}.  Under these conditions,
we describe the various possible $A_{\infty}$-structures on the Ext
algebras.  In Section~\ref{yysect6}, we eliminate some of those
possibilities, showing that they correspond to non-regular algebras.
In Section~\ref{yysect-reg}, we show that the remaining possibilities
lead to regular algebras, and this leads to the proofs of our main
theorems, which we give in Section~\ref{yysect7}.

There is also an appendix, which gives results on
$A_\infty$-structure on Ext-algebras.  These have appeared
elsewhere, some without proof.  The appendix starts with a discussion of
Kadeishvili's and Merkulov's results about the
$A_{\infty}$-structure on the homology of a DGA.  The bar construction
is described: this is a DGA whose homology is Ext, and so leads to an
$A_{\infty}$-structure on an Ext algebra.  The appendix finishes with
a proof of Keller's theorem.  Throughout, some extra care is taken
with gradings, since this is necessary in various parts of this paper.

\section{Artin-Schelter regular algebras}
\label{yysect1}

In this section, we recall the definition of Artin-Schelter regular
algebras, we study their Hilbert series, and for those of global
dimension 4 generated in degree 1, we describe the possible shapes of
their Ext-algebras.

A connected graded algebra $A$ is called {\it Artin-Schelter regular} 
(or {\it AS regular}) if the following three conditions hold.
\begin{enumerate}
\item[(AS1)] $A$ has finite global dimension $d$,
\item[(AS2)] $A$ is {\it Gorenstein}, i.e., for some integer $l$,
\[
\Ext^i_A(k, A)=\begin{cases} k(l) & \text{ if } i=d\\
				0   & \text{ if } i\neq d
\end{cases}
\]
where $k$ is the trivial module $A/\fm$, and
\item[(AS3)] 
$A$ has finite Gelfand-Kirillov dimension, i.e., there is a 
positive number $c$ such that $\dim A_n< c\; n^c$ for all 
$n\in \mathbb{N}$.
\end{enumerate}

The notation $(l)$ in (AS2) is the degree $l$ shift operation on
graded modules.

AS regular algebras have been studied in many recent papers, 
and in particular, AS regular algebras of global dimension 3 have 
been classified \cite{ASc, ATV1, ATV2, Ste1, Ste2} and their geometry 
has been studied extensively  \cite{ATV2, Ste3}. If $A$ is an AS 
regular algebra of global dimension 3, then it is generated by two 
or three elements. Suppose now that $A$ is generated in degree 1. If $A$ 
is (minimally) generated by three elements, then $A$ is Koszul and the 
trivial $A$-module $k$ has a minimal free resolution  of the form
\[
0\to A(-3)\to A(-2)^{\oplus 3}\to A(-1)^{\oplus 3}\to A\to k\to 0.
\]
If $A$ is generated by two elements, then $A$ is not Koszul and the trivial
$A$-module $k$ has a minimal free resolution of the form
\[
0\to A(-4)\to A(-3)^{\oplus 2}\to A(-1)^{\oplus 2}\to A\to k\to 0.
\]
In this case, since $A$ is not Koszul, the Ext-algebra of $A$ has
nontrivial higher multiplications. The Ext-algebra of $A$ has been
studied by Shi and Wang \cite{SW}.


\begin{lemma}
\label{yylem1.1}
Suppose $A$ is connected graded and satisfies \textup{(AS1)} and 
\textup{(AS2)}.
\begin{enumerate}
\item $A$ is finitely generated.
\item The trivial $A$-module $k_A$ has a minimal free resolution of the form
\[
0\to P_{d}\to \cdots \to P_{1}\to P_{0}\to k_A\to 0,
\]
where $P_{w}=\boplus_{s=1}^{n_w}A(-i_{w,s})$ for some finite integers 
$n_w$ and $i_{w,s}$.
\item The above free resolution is symmetric in the following sense:
$P_0=A$, $P_{d}=A(-l)$, $n_w=n_{d-w}$, and $i_{w,s}+i_{d-w, n_w-s+1}=l$
for all $w,s$.
\end{enumerate}
\end{lemma}

\begin{proof} (a) and (b) are proved in \cite[3.1.1]{SZ}.
Part (c) is given in the proof of \cite[3.1.4]{SZ}. We 
repeat the main idea here. The trivial $A$-module $k_A$ has a minimal free 
resolution of  the form
\begin{equation}
\label{E1.1.1}
0\to P_{d}\to \cdots \to P_{1}\to P_{0}\to k_A\to 0,
\end{equation}
where $P_{w}=\boplus_{s=1}^{n_w}A(-i_{w,s})$ for some finite integers 
$n_w$ and $i_{w,s}$. The dual complex $(\textup{\ref{E1.1.1}})^{\vee}$ 
is a free  resolution of $_Ak(l)$. By duality, \eqref{E1.1.1} is minimal 
if and only if  $(\textup{\ref{E1.1.1}})^{\vee}$ is. Since the $i$th 
term of the minimal free resolution of $k_A$ (respectively, $_Ak$) is 
isomorphic to  $\Tor^A_i(k,k)\otimes_k A$ (respectively, 
$A\otimes_k \Tor^A_i(k,k)$), $k_A$ and $_Ak$ have the same type of 
minimal free resolution.  By comparing \eqref{E1.1.1} with 
$(\textup{\ref{E1.1.1}})^{\vee}(-l)$,  we see that $-i_{w,s}
=-l+i_{d-w,n_w-s+1}$. It is clear that $P_0=A$, so $P_d=A(-l)$.
\end{proof}

Recall from \cite{AZ} that the cohomological dimension of the 
noncommutative projective scheme $\Proj A$ associated to a graded ring $A$
is defined to be
\[
\cd \Proj A=\sup\{i\;|\; \HB^i(\Proj A,\mathcal{M})\neq 0 \ 
\textrm{for some} \  \mathcal{M}\in \QGr A\}.
\]

\begin{lemma}
\label{yylem1.2} If $A$ is a noetherian AS regular algebra of 
global dimension at least 3, then the GK-dimension of $A$ is at 
least $3$.
\end{lemma}

\begin{proof} Let $A$ be a noetherian connected graded ring with
$\GKdim A\leq 2$. We claim that $\cd \Proj A\leq 1$. Since $A$ is 
noetherian, we may assume that $A$ is prime by 
\cite[8.5]{AZ}. If $\GKdim A\leq 1$, then by \cite{SSW}, 
$A$ is PI and hence by \cite[8.13]{AZ} $\cd \Proj A\leq 0$. 
So it suffices to consider the case when $\GKdim A=2$. Combining the 
results of Artin-Stafford \cite[0.3, 0.4 and 0.5]{ASt2}, we 
obtain that some Veronese subring $A^{(n)}$ of $A$ is a subring of a 
twisted homogeneous coordinate ring $B:=B(\mathcal{E},\mathcal{B}_1,\tau)$ 
and $\Proj A^{(n)}\cong \Proj B$. Here $\mathcal{E}$ is an 
$\mathcal{O}_Y$-order for some projective curve $Y$ and $\mathcal{B}_1$ is 
an ample invertible $(\mathcal{E},\mathcal{E}^\tau)$-bimodule; see 
\cite{ASt2} for the details. In particular, \cite[0.4(ii)]{ASt2}
says that $\Proj A^{(n)}$ and $\Mod \mathcal{O}_{\mathcal E}$ are
equivalent. By \cite[8.7(2)]{AZ}, $\Proj A=\Proj A^{(n)}$. 
Hence 
$$\cd \Proj A=\cd \Proj A^{(n)}=\cd \Mod \mathcal{O}_{\mathcal E}=1,$$ 
where the last equality follows from the fact that $Y$ is a curve.

On the other hand, by \cite[8.1]{AZ}, if $A$ is AS regular, 
then $\cd \Proj A=d-1$ where $d$ is the global dimension of $A$. This 
shows that if $d\geq 3$, then $\cd \Proj A\geq 2$. By the 
claim proved in the last paragraph one sees that $\GKdim A>2$. 
Finally, the GK-dimension of a noetherian AS regular algebra is an 
integer \cite[2.4]{SZ}, so $\GKdim\geq 3$.
\end{proof}

Now we start to focus on AS regular algebras of global 
dimension 4.  The classification of such algebras is far from 
finished.  Some abstract properties of 4-dimensional AS 
regular algebras have been proved.  We will use the following result
from \cite[3.9]{ATV2} in later sections.

\begin{lemma}
\label{yylem1.3}
\cite[3.9]{ATV2}.
If $A$ is a noetherian connected graded AS regular algebra of
global dimension 4, then it is an integral domain.
\end{lemma}

The Hilbert series of a graded vector space $M=\boplus_{i\in {\mathbb Z}}
M_i$ is defined to be 
\[
H_M(t)=\sum_{i\in {\mathbb Z}}(\dim_k M_i) t^i.
\]
We determine the Hilbert series of $A$ when $A$ is generated in degree 1. 

\begin{proposition} 
\label{yyprop1.4}
Let $A$ be a graded AS regular algebra of global dimension 4 
that is generated in degree 1.  Suppose that $A$ is a domain.
Then $A$ is minimally generated by either 2, 3, or 4 elements. 
\begin{enumerate}
\item
If $A$ is generated by 2 elements, then 
there are two relations whose degrees are 3 and 4. 
The minimal resolution of the trivial module is of the form
\[
0\to A(-7)\to A(-6)^{\oplus 2}\to A(-4)\oplus A(-3)\to A(-1)^{\oplus 2}
\to A\to k\to 0.
\]
The
Hilbert series of $A$ is 
\[
H_A(t)={1/ {(1-t)^2 (1-t^2)(1-t^3)}}.
\]
\item 
If $A$ is generated by 3 elements, then 
there are two relations in degree 2 and two relations in degree 3. 
The minimal resolution of the trivial module is of the form
\[
0\to A(-5)\to A(-4)^{\oplus 3}\to A(-3)^{\oplus 2}\oplus A(-2)^{\oplus 2}
\to A(-1)^{\oplus 3}
\to A\to k\to 0.
\]
The Hilbert series of $A$ is 
\[
H_A(t)={1/ (1-t)^3 (1-t^2)}.
\]
\item
If $A$ is generated by 4 elements, then there are six quadratic relations.
The minimal resolution of the trivial module is of the form
\[
0\to A(-4)\to A(-3)^{\oplus 4}\to A(-2)^{\oplus 6}\to A(-1)^{\oplus 4}
\to A\to k\to 0.
\]
The Hilbert series of $A$ is 
\[
H_A(t)={1/ (1-t)^4}.
\]
\end{enumerate}
In each of these cases, the GK-dimension of $A$ is $4$.
\end{proposition}

\begin{proof} Suppose $A$ is minimally generated by $n$ elements in 
degree 1. By Lemma \ref{yylem1.1}, $k_A$ has a minimal free resolution
\[
0\to A(-l)\to A(-l+1)^{\oplus n}\to \boplus_{s=1}^v A(-n_s)\to 
A(-1)^{\oplus n} \to A\to k\to 0 
\]
where $2\leq n_1\leq \cdots \leq n_v\leq l-2$ and $n_s+n_{v-s+1}=l$ for all 
$s$. The Hilbert series of $A$ is $H_A(t)=1/p(t)$, where
\[
p(t)=1-nt+\sum_{s=1}^v t^{n_s}-nt^{l-1}+t^l.
\]
Recall that the GK-dimension of $A$ is equal to the order of the zero 
of $p(t)$ at $1$ \cite[2.2]{SZ}. Clearly $\GKdim A\geq 1$, whence 
$p(1)=0$. This  implies that $v=2n-2$.
Since $n_s+n_{v-s+1}=l$, we have 
\begin{equation}
\label{E1.4.1}
\sum_s n_s={\frac{1}{2}}\sum_{s=1}^v l={\frac{1}{2}}(2n-2)l=(n-1)l.
\end{equation}
We claim that $\GKdim A\geq 3$.  If $\GKdim A=1$, then $A$ is PI and 
noetherian by \cite{SSW}, but if $A$ is noetherian, then by
Lemma~\ref{yylem1.2}, $\GKdim A \geq 3$.  If $\GKdim A=2$, by
\cite{ASt1} $A$ is isomorphic to a twisted homogeneous coordinate ring 
$B:=B(E,\sigma, {\mathcal L})$ up to a finite-dimensional vector 
space, where $E$ is a projective curve. In particular, $A$ is 
noetherian, and again Lemma~\ref{yylem1.2} yields a contradiction.
Therefore $\GKdim A\geq 3$. Now going back to 
the polynomial $p(t)$, the fact that $\GKdim A\geq 3$ implies that $p''(1)=0$.
Since
\[
p''(t)=\sum_{s=1}^v n_s(n_s -1) t^{n_s-2}-
(l-1)(l-2)n t^{l-3}+l(l-1)t^{l-2},
\]
we obtain
\[
\sum_{s=1}^v n_s(n_s -1)-n(l-1)(l-2)+l(l-1)=0.
\]
Simplifying this by using \eqref{E1.4.1}, we have
\begin{equation}
\label{E1.4.2}
\sum_{s=1}^v n_s^2=n(l^2-2l+2)-l^2.
\end{equation}
Using \eqref{E1.4.1} and \eqref{E1.4.2}, we have
\[
\sum_{s=1}^v n_{v-s+1}n_s=\sum_s (l-n_s)n_s=l\sum_s n_s-\sum_s n_s^2
=2(l-1)n.
\]
Since $v=2(n-1)$, we have
\begin{equation}
\label{E1.4.3}
\sum_{s=1}^{n-1} n_{v-s+1}n_s=(l-1)n.
\end{equation}
We will do some detailed analysis using \eqref{E1.4.1}---\eqref{E1.4.3}.

Clearly $n\geq 2$. If $n=2$, then $v=2$, and \eqref{E1.4.3} becomes $n_1(l-n_1)
=2(l-1)$. We see that $n_1\neq 2$. Hence $n_1\geq 3$, and since
$n_1\leq l/2$, $n_1(l-n_1)\geq 3(l-3)$. Then $3(l-3)\leq 2(l-1)$ implies 
that $l\leq 7$. Since $n_1\leq l/2$, $n_1=3$. So $n_2=3$ (if $l=6$) or
$n_2=4$ (if $l=7$). 
Also $l=n_1+n_2\geq 6$. Since $n_1 n_2\neq 2\cdot 5$, $l\neq 6$. Hence
$l=7$ is the only possible value.  When $l=7$, then $n_1=3$ and $n_2=4$.
(a) follows. 

Now we consider the case when $n\geq 3$. Since $2\leq n_s\leq {\frac{1}{2}l}$ 
for $1\leq s\leq n-1$, we have
\[
2(l-2)\leq n_s(l-n_s),
\]
and \eqref{E1.4.3} implies that
\[
(n-1)2(l-2)\leq  n(l-1).
\]
Solving the inequality, we have
\[
l\leq \frac{3n-4}{n-2}=3+ \frac{2}{n-2}.
\]

If $n=3$, then the last inequality implies that $l\leq 5$. Also 
$l\geq 4$ for  any AS regular algebra of global dimension 4. If $l=4$ 
(and $n=3$), by \eqref{E1.4.1} 
and the fact that $n_s\geq 2$, we find that $n_s=2$ for all $s$.  This contradicts 
\eqref{E1.4.3}, though. Hence $l=5$ and $n_1=n_2=2$, $n_3=n_4=3$. (b) follows.

If $n=4$, $l\leq 4$ and hence $l=4$. Thus $n_s=2$ for all $s$, $1 \leq
s \leq 6$.  (c) follows.

If $n\geq 5$, then 
\[
l\leq 3+ \frac{2}{n-2}<4.
\] 
This is impossible for an AS regular algebra of global dimension 4. 

It is clear from the Hilbert series that $\GKdim A=4$ in all
three cases.
\end{proof}

Artin and Schelter stated the same result in the PI case \cite[1.20]{ASc}
and provided some examples in each of the cases (a,b,c).

\begin{remark} 
\label{yyrem1.5} In Proposition \ref{yyprop1.4} it is essential to
assume that $A$ has finite GK-dimension, i.e., (AS3). 
If $A$ does not have finite GK-dimension, then $A$ has
exponential growth and the number of the generators can be any number 
larger than 4.  This means that there are many more algebras satisfying 
the other hypotheses of Proposition \ref{yyprop1.4}.
\end{remark}

We do not have a counterexample to the following question.

\begin{question}
\label{yyque1.6} 
If $A$ is noetherian and connected graded, is the minimal 
number of generators of $A$ no more than the global dimension 
of $A$?
\end{question}

\begin{definition}
\label{yydefn1.7}
If an algebra satisfies the hypotheses of Proposition~\ref{yyprop1.4},
we label it according to the dimensions of vector spaces
$\Ext^i_A(k,k)$.  That is, algebras as in Proposition
\ref{yyprop1.4}(a) are said to be of \emph{type \type{12221}},
algebras as in Proposition \ref{yyprop1.4}(b) are of {\it type \type{13431}},
and algebras as in Proposition \ref{yyprop1.4}(c) are of {\it type \type{14641}}.
\end{definition}


We also state the following result for later use.  This is a
generalization of a result of Smith \cite{Sm1}, and was proved using
$A_\infty$-algebra methods.

\begin{theorem} 
\label{yythm1.8} \cite{LPWZ1} 
Let $A$ be a connected graded algebra and let $E$ be the Ext-algebra of $A$. 
Then $A$ is AS regular if and only if $E$ is Frobenius.
\end{theorem}

\section{$A_\infty$-algebras and grading} 
\label{yysect2}

The notion of $A_\infty$-algebra was introduced by Stasheff in the
1960s \cite{Sta}.  Since then, more and more theories involving
$A_\infty$-structures have been discovered in algebra, geometry and
mathematical physics; see Kontsevich's 1994 ICM paper \cite{Ko}, for
example.  The methods of $A_\infty$-algebras can be very effective in
ring theory. The best introductory paper on $A_\infty$-algebras for
ring theorists is a paper by Keller \cite{Ke3}; Keller's other papers
\cite{Ke2,Ke4,Ke5} are also very interesting. These papers convinced
us to look further into applications of $A_\infty$-algebras in ring
theory.

In this section we review the definition of an $A_{\infty}$-algebra,
we discuss grading systems, and we state a result about
$A_{\infty}$-structures on Ext algebras.  Other basic material 
about $A_\infty$-algebras can be found in Keller's paper \cite{Ke3}.
Some examples of $A_\infty$-algebras related to ring theory were given in 
\cite{LPWZ2}. 

\subsection{$A_{\infty}$-algebras}

\begin{definition}
\label{yydef2.1} \cite{Sta}
An {\it $A_\infty$-algebra} over a base field $k$ is a $\mathbb Z$-graded 
vector space
\[
A=\boplus_{p\in {\mathbb Z}} A^p
\]
endowed with a family of graded $k$-linear maps
\[
m_n: A^{\otimes n} \to A, \quad n\geq 1,
\]
of degree $2-n$ satisfying the following {\it Stasheff identities}:
\begin{equation}\label{SI}
\sum (-1)^{r+st} m_u(id^{\otimes r}\otimes m_s \otimes id^{\otimes
t})=0
\tag*{\SI{n}}
\end{equation}
for all $n\geq 1$, where the sum runs over all decompositions 
$n=r+s+t$ ($r,t\geq 0$ and $s\geq 1$), and where $u=r+1+t$. Here, 
$id$ denotes the identity map of $A$. Note that when these formulas 
are applied to elements, additional  signs appear due to the Koszul 
sign rule. Some authors also use the terminology {\it strongly 
homotopy associative algebra} (or {\it sha algebra}) for 
$A_\infty$-algebra.
\end{definition}

The degree of $m_1$ is $1$ and the identity \SI{1} is $m_1 m_1=0$. This
says that $m_1$ is a differential of $A$. The identity \SI{2} is
\[
m_1 m_2= m_2(m_1\otimes id+id\otimes m_1)
\]
as maps $A^{\otimes 2} \to A$. So the differential
$m_1$ is a graded derivation with respect to $m_2$. Note that $m_2$
plays the role of multiplication although it may not be associative.
The degree of $m_2$ is zero. The identity \SI{3} is
\[
m_2(id\otimes m_2-m_2\otimes id)
=m_1m_3+ m_3(m_1\otimes id \otimes id+ id\otimes m_1\otimes
id + id\otimes id \otimes m_1)
\]
as maps $A^{\otimes 3}\to A$. If either $m_1$ or $m_3$ is zero, then
$m_2$ is associative. In general, $m_2$ is associative up to a
chain homotopy given by $m_3$.

When $n\geq 3$, the map $m_n$ is called a {\it higher multiplication}. We 
write  an $A_\infty$-algebra $A$ as $(A, m_1, m_2, m_3, \cdots)$ to indicate
the multiplications $m_i$. An associative algebra $A$ is an 
$A_\infty$-algebra concentrated in degree 0 with all multiplications $m_n=0$ 
for $n\neq 2$, so an associative algebra has the form of $(A,m_2)$.  A 
differential graded algebra (or DGA) $(A,\partial)$ is an 
$A_\infty$-algebra with $m_1=\partial$, $m_2$ the multiplication,
and $m_n=0$ for all $n\geq 3$, so we may write a DGA as $(A,m_1, m_2)$. 
We also assume that every $A_\infty$-algebra in this paper contains an 
identity element $1$ with respect to the multiplication $m_2$ that satisfies 
the following {\it strictly unital condition}:

\bigskip

If $n\neq 2$ and $a_i=1$ for some $i$, then $m_n(a_1,\cdots, a_n)=0$. 

\bigskip
\noindent
In this case, $1$ is called the {\it strict unit} or {\it identity} of $A$.

Let $A$ be a graded algebra generated in degree 1, and let $k_A$ be the trivial 
$A$-module. Then the Ext-algebra $\Ext^*_A(k_A,k_A)$ is equipped with 
an $A_\infty$-algebra structure. We use 
$\Ext^*_A(k_A,k_A)$ to denote both the usual associative Ext-algebra
and the Ext-algebra with its $A_\infty$-structure.
For more information on this $A_{\infty}$-structure, see
Subsection~\ref{sectapp3}.  By \cite[Ex. 13.4]{LPWZ2} there is a
graded algebra $A$ such that the associative algebra
$\Ext^*_A(k_A,k_A)$ does not contain enough information to recover the
original algebra $A$; on the other hand, the information from the
$A_\infty$-algebra $\Ext^*_A(k_A,k_A)$ is sufficient to recover $A$.
This is the point of the following theorem, and this process of
recovering the algebra from its Ext-algebra is one of the main tools
used in this paper.

If $V$ is a graded vector space over $k$, then the graded $k$-linear
dual of $V$ is denoted by $V^{\#}$.

\begin{theorem}
[Keller's higher-multiplication theorem in the connected 
graded case] 
\label{yythm2.2}
Let $A$ be a graded algebra, finitely generated in degree 1, and let $E$ be the 
$A_\infty$-algebra $\Ext^*_A(k_A,k_A)$. Let $R=\boplus_{n\geq 2} R_n$ be the 
minimal graded space of relations of $A$ such that $R_n\subset A_1\otimes 
A_{n-1}\subset A_1^{\otimes n} $. Let $i: R_n\to A_1^{\otimes n}$ be the 
inclusion map and let $i^{\#}$ be its $k$-linear dual. Then the 
multiplication $m_n$ of $E$ restricted to $(E^1)^{\otimes n}$ is equal 
to the map
\[
i^{\#}: (E^1)^{\otimes n}=(A_1^{\#})^{\otimes n}\longrightarrow 
R^{\#}_n\subset E^2.
\]
\end{theorem}

Keller has the same result for a more general class of algebras $kQ/I$ 
where $Q$ is a finite quiver and $I$ is an admissible ideal of $kQ$;
this was stated in \cite[Proposition 2]{Ke4} without proof.  Our result 
works only for graded algebras generated in degree 1, so it is a special 
case of Keller's result.  We announced Theorem \ref{yythm2.2} at several 
conferences a few years ago, where several experts informed us that the
result was somewhat known and they were interested in a detailed proof.
As we stated in \cite{LPWZ2}, we hope to develop a theory of 
$A_\infty$-algebras for ring theorists, especially for people working 
on homological properties of graded rings.  Keeping this in mind, a 
quite detailed proof of Theorem \ref{yythm2.2} is provided in the Appendix. 
We also want to point out that Theorem \ref{yythm2.2} is essential for 
the classification of Artin-Schelter regular algebras of global 
dimension 4 that are generated by two elements \cite{LWZ}. 

We mention that there may be several quasi-isomorphic
$A_{\infty}$-algebra structures on $\Ext_{A}^{*}(k_{A}, k_{A})$; we call
these different structures \emph{models} for the quasi-isomorphism
class of $E$.

\subsection{Adams grading}

In this paper we are mainly interested in graded algebras and their
Ext-algebras.  The grading appearing in a graded algebra may be
different from the grading appearing in the definition of the
$A_\infty$-algebra.  We introduce the Adams grading for an
$A_\infty$-algebra, as follows. Let $G$ be an abelian group.  (In this
paper, $G$ will always be free abelian of finite rank.)  Consider a
bigraded vector space
\[
A=\boplus_{p\in \Z, i\in G} A^p_i
\]
where the upper index $p$ is the grading appearing in Definition
\ref{yydef2.1}, and the lower index $i$ is an extra grading, called
the \emph{$G$-Adams grading}, or {\it Adams grading} if $G$ is
understood.  We also write
\[
A^p=\boplus_{i\in G}A^p_i \quad\text{and}\quad A_i=
\boplus_{p\in \Z}A^p_i.
\]
The degree of a nonzero element in $A^p_i$ is $(p,i)$, and the second
degree is called the {\it Adams degree}.  For an $A_{\infty}$-algebra
$A$ to have an Adams grading, the map $m_n$ in Definition
\ref{yydef2.1} must be of degree $(2-n,0)$: each $m_n$ must preserve
the Adams grading. When $A$ is an associative $G$-graded algebra
$A=\boplus_{i\in G} A_i$, we view $A$ as an $A_\infty$-algebra (or a
DGA) concentrated in degree 0, viewing the given grading on $A$ as the
Adams grading.  The Ext-algebra of a graded algebra is bigraded; the
grading inherited from the graded algebra is the Adams grading, and we
keep using the lower index to denote the Adams degree.

Assume now that $G=\Z$, since we are mainly interested in this case.
Write 
\[
A^{\geq n}=\boplus_{p\geq n}A^p \quad\text{and}\quad A_{\geq n}=
\boplus_{i\geq n}A_i,
\]
and similarly for $A^{\leq n}$ and $A_{\leq n}$.
An
$A_\infty$-algebra $A$ with a $\Z$-Adams grading
is called {\it Adams connected} if (a) $A_0=k$, (b) $A=A_{\geq 0}$ 
or $A=A_{\leq 0}$, and (c) $A_i$ is finite-dimensional for all $i$. 
When $G=\Z\times G_0$, we define {\it Adams connected} 
in the same way after omitting the $G_0$-grading. If $A$ is a connected 
graded algebra which is finite-dimensional in each degree, then it is 
Adams connected when viewed as an $A_\infty$-algebra concentrated in 
degree 0. 

We prove the following result in the appendix.

\begin{proposition}
\label{yyprop2.3} Let $A$ be a $\Z\oplus G$-Adams graded algebra,
such that with 
respect to the $\Z$-grading, $A$ is finitely generated in degree 1. 
Then there is an $A_{\infty}$-model for $E$ such that the multiplication $m_n$ 
in Theorem~\ref{yythm2.2} preserves the $\Z\oplus G$ grading.
\end{proposition}

\section{$A_\infty$-Ext-algebras of regular algebras of dimension 4}
\label{yysect3}

Let $E$ be the $A_\infty$-Ext-algebra $\Ext^*_A(k,k)$ of $A$.
For simplicity, we just call it the $\Ext$-algebra of $A$.  In this
section we use information about the grading of $E$ from
Proposition~\ref{yyprop1.4} to examine the possible
$A_{\infty}$-structures on $E$.

\begin{proposition} 
\label{yyprop3.1}
Let $A$ be an algebra as in \textup{Proposition \ref{yyprop1.4}} and 
let $E$ be the $\Ext$-algebra of $A$.
\begin{enumerate}
\item
(type \type{12221})\ \ 
If $A$ is minimally generated by 2 elements, then $E$ is isomorphic to 
\[
k\boplus E^{1}_{-1}\boplus (E^{2}_{-3}\oplus E^{2}_{-4}) 
\boplus E^{3}_{-6}\boplus E^{4}_{-7}
\]
as a ${\Z}^2$-graded vector space, where the lower index is the
Adams grading inherited from the grading of $A$ and the upper index is 
the homological grading of the Ext-group.  The dimensions of the subspaces are
\[
\dim E^1_{-1}=\dim E^{3}_{-6}=2, \qquad \dim E^{2}_{-3}=\dim E^{2}_{-4}
=\dim E^{4}_{-7}=1.
\]
As an $A_\infty$-algebra, $m_n=0$ for all
$n\geq 5$; that is, $E=(E,m_2,m_3,m_4)$.
\item
(type \type{13431})\ \ 
If $A$ is minimally generated by 3 elements, then $E$ is isomorphic to 
\[
k\boplus E^{1}_{-1}\boplus (E^{2}_{-2}\oplus E^{2}_{-3}) 
\boplus E^{3}_{-4}\boplus E^{4}_{-5}
\]
as a ${\Z}^2$-graded vector space. As an $A_\infty$-algebra, 
$m_n=0$ for all $n\geq 4$; that is, $E=(E,m_2,m_3)$.
The dimensions of the subspaces are
\[
\dim E^1_{-1}=\dim E^{3}_{-4}=3, \quad \dim E^{2}_{-2}=\dim E^{2}_{-3}=2,
\quad \dim E^{4}_{-5}=1.
\]
\item
(type \type{14641})\ \ 
If $A$ is minimally generated by 4 elements, then $E$ is isomorphic to 
\[
k\boplus E^{1}_{-1}\boplus E^{2}_{-2} \boplus E^{3}_{-3}\boplus 
E^{4}_{-4}
\]
as a ${\Z}^2$-graded vector space.  The algebras $A$ and $E$ are Koszul
and $m_n$ of $E$ is zero for all $n\neq 2$.  The dimensions of the subspaces are
\[
\dim E^1_{-1}=\dim E^{3}_{-3}=4, \qquad \dim E^{2}_{-2}=6,
\quad \dim E^{4}_{-4}=1.
\]
\end{enumerate}
\end{proposition}

\begin{proof} The vector space decomposition of $E$ and the dimensions of 
the subspaces of $E$ are clear from the form of the minimal free resolution 
of the trivial module. In case (c), $k_A$ has a linear resolution
in the sense of \cite{Sm1}. Hence $A$ 
and $E$ are Koszul, and this implies that the higher multiplications
of $E$ are trivial.  The assertions about higher multiplications
in (a) and (b) (and also (c)) follow from the next lemma.
\end{proof}

\begin{lemma} 
\label{yylem3.2} Let $E$ be an Adams connected $A_\infty$-algebra. 
\begin{enumerate}
\item
(type \type{12221})\ \ 
If $E$ is isomorphic to 
\[
k\boplus E^{1}_{-1}\boplus (E^{2}_{-3}\oplus E^{2}_{-4}) 
\boplus E^{3}_{-6}\boplus E^{4}_{-7}
\]
as a ${\Z}^2$-graded vector space, then 
$m_n=0$ for $n\neq 2,3,4$. 
\item
(type \type{13431})\ \ 
If $E$ is isomorphic to 
\[
k\boplus E^{1}_{-1}\boplus (E^{2}_{-2}\oplus E^{2}_{-3}) 
\boplus E^{3}_{-4}\boplus E^{4}_{-5}
\]
as a ${\Z}^2$-graded vector space, then $m_n=0$ for $n\neq 2,3$.
\item
(type \type{14641})\ \ 
If $E$ is isomorphic to 
\[
k\boplus E^{1}_{-1}\boplus E^{2}_{-2} \boplus E^{3}_{-3}\boplus 
E^{4}_{-4}
\] 
as a ${\Z}^2$-graded vector space, then $m_n=0$ for $n\neq 2$.
\end{enumerate}
\end{lemma}

\begin{proof} First of all $m_1=0$ since the degree of $m_1$ is $(1,0)$.

(a) We consider the maps $m_n$ for $n \geq 5$, restricted to 
a homogeneous subspace:
\[
m_n: E^{i_1}_{-j_1}\otimes \cdots \otimes E^{i_n}_{-j_n}\to E^i_{-j}
\]
where $i=\sum_{s=1}^n i_s -n+2$ and $j=\sum_{s=1}^n j_s$. If $i_s=0$ 
(and hence $E^{i_s}_{-j_s}=k$) for some $s$, then $m_n=0$ by the strict
unital condition. So we assume that $i_s\geq 1$ for all $s$. 

If $E^{i_s}_{-j_s} \neq 0$, then $i_s\leq j_s\leq i_s +3$. We assume that 
$E^{i_s}_{-j_s}\neq 0$ for all $s$ and show that $E^{i}_{-j}=0$ when 
$n\geq 5$. If $n>5$, then 
\[
j=\sum_{s} j_s\geq \sum_s i_s =i+(n-2)> i+3.
\]
Hence $E^{i}_{-j}=0$ and $m_n=0$. It remains to show that $m_5=0$. If $i_s=1$
for all $s$, then $E^{i}_{-j}=E^2_{-5}=0$, whence $m_5=0$. If $i_{s_0}>1$ 
for some $s_0$, then $j_{s_0}\geq i_{s_0} +1$ and 
\[
j=\sum_{s} j_s\geq \sum_s i_s +1 =i+(5-2)+1> i+3.
\]
This shows that $E^i_{-j}=0$ and $m_5=0$.

(b) The proof is similar to that of (a). We again consider the 
multiplications $m_n$ restricted to a homogeneous subspace:
\[
m_n: E^{i_1}_{-j_1}\otimes \cdots \otimes E^{i_n}_{-j_n}\to E^i_{-j}.
\]
We assume that $E^{i_s}_{-j_s}\neq 0$ for all $s$ and show that 
$E^{i}_{-j}=0$ when $n\geq 4$. 

If $E^{i_s}_{-j_s} \neq 0$, then $i_s\leq j_s< i_s +2$. If $n\geq 4$, then 
\[
j=\sum_{s} j_s\geq \sum_s i_s =i+(n-2)\geq i+2.
\]
Hence $E^{i}_{-j}=0$ and $m_n=0$. 

(c) Use an argument similar to (a) or (b).
\end{proof}

If $A$ is an AS regular algebra of type \type{14641}, then $A$ and its 
Ext-algebra are both Koszul. All higher multiplications of
$E$ are zero, so $A_\infty$-algebra methods cannot be used here.
We will have other methods to study these algebras, for example, 
using the matrices constructed from the multiplication of $E$. 
We present this idea below for AS regular algebras of type 
\type{13431}.

If $A$ is an AS regular algebra of type \type{13431}, then $m_2$ and $m_3$ 
of $E$ will be nonzero. By Theorem \ref{yythm2.2} and Proposition
\ref{yyprop3.1}, the relations of $A$, two of which are in degree 2
and other two are in degree 3, are determined completely 
by $m_2$ and $m_3$ of $E$.  By the ${\mathbb Z}^2$-graded
decomposition given in Lemma \ref{yylem3.2}(b), the possible nontrivial
components of $m_2$ of $E$ are
\begin{gather*}
E_{-1}^1\otimes E_{-1}^1\to E_{-2}^2,  \\
E_{-1}^1\otimes E_{-3}^2\to E_{-4}^3, 
\quad E_{-3}^2\otimes E_{-1}^1\to E_{-4}^3, \\
E_{-1}^1\otimes E_{-4}^3\to E_{-5}^4, 
\quad E_{-4}^3\otimes E_{-1}^1\to E_{-5}^4, \\
E_{-2}^2\otimes E_{-3}^2\to E_{-5}^4, 
\quad E_{-3}^2\otimes E_{-2}^2\to E_{-5}^4.
\end{gather*}
Let $\{\alpha_1,\alpha_2,\alpha_3\}$ be a $k$-linear basis of 
$E_{-1}^1$. Let $\{\beta_1,\beta_2\}$ be a $k$-linear basis of
$E_{-2}^2$ and $\{\zeta_1,\zeta_2\}$ a $k$-linear basis of
$E_{-3}^2$. Let $\{\gamma_1,\gamma_2,\gamma_3\}$ be a $k$-linear
basis of $E_{-4}^3$, and let $\delta$ be a nonzero element of
$E_{-5}^4$.  Since $E$ is Frobenius (see Theorem \ref{yythm1.8}),
the maps
\[
m_2: E_{-1}^1\otimes E_{-4}^3\to E_{-5}^4 \quad \text{and} \quad 
m_2: E_{-4}^3\otimes E_{-1}^1\to E_{-5}^4
\] 
define a perfect pairing. We may choose $\{\gamma_1,\gamma_2,\gamma_3\}$
so that $\alpha_i\gamma_j=\delta_{ij}\delta$.
Write $\gamma_i\alpha_j=\lambda_{ij}\delta$; then $\Lambda:
=(\lambda_{ij})_{3\times 3}$ is a non-singular matrix.
Similar we may choose $\{\zeta_1,\zeta_2\}$ such that $\beta_i\zeta_j=
\delta_{ij}\delta$. Write $\zeta_i\beta_j=t_{ij}\delta$; then
the matrix $T:=(t_{ij})_{2\times 2}$ is a non-singular matrix.
We call $(\Lambda,T)$ the {\it Frobenius data} of $E$. 
For algebras of type \type{14641}, the Frobenius data provides a lot of
information about $A$.  For algebras of type \type{13431}, one uses the
different Jordan forms of $(\Lambda,T)$ and the associativity of $m_2$
to obtain a list of equations which the entries of $\Lambda$ and $T$
must satisfy.  These equations should be essential for the
classification of the possible $m_2$'s.

Next we look into $m_3$. Here is a 
list of the possible nonzero components of $m_3$:
\begin{gather*}
(E_{-1}^1)^{\otimes 3}\to E_{-3}^2, \\
[(E_{-1}^1)^{\otimes 2}\otimes E_{-2}^2]\oplus
[E_{-1}^1\otimes E_{-2}^2\otimes E_{-1}^1]\oplus
[E_{-2}^2\otimes (E_{-1}^1)^{\otimes 2}]\to E_{-4}^3, \\
E_{-1}^1\otimes (E_{-2}^2)^{\otimes 2}]\oplus
[E_{-2}^2\otimes E_{-1}^1\otimes E_{-2}^2]\oplus
[(E_{-2}^2)^{\otimes 2}\otimes E_{-1}^1]\to E_{-5}^4.
\end{gather*}
It is not hard to verify that the identities \SI{n} are automatic for
$n\geq 5$. The only constraint to make $(E,m_2,m_3)$ an
$A_\infty$-algebra is \SI{4}. In the near future we will look into the
possibility of classifying those $A_{\infty}$-algebras $(E,m_2,m_3)$
which correspond to the AS regular algebras of type \type{13431}.

\section{$A_\infty$-Ext-algebras of type \type{12221}}
\label{yysect4}

We now concentrate on algebras of type \type{12221}.  In this section
we describe formulas for the possible multiplication maps $m_{n}$ on
their $A_{\infty}$-Ext-algebras, and we use Keller's
theorem~\ref{yythm2.2} to relate these formulas to the relations in
the original algebra.  This sets the stage for our classification,
which starts in the next section.

Using the ideas in the proof of Lemma \ref{yylem3.2}, we can list all
of the possible maps $m_n$ here. 
Except for the multiplying by the unit element, the possible nonzero $m_2$'s of $E$ are
\begin{gather*}
E^{1}_{-1}\otimes E^{3}_{-6}\to E^{4}_{-7}, \qquad 
E^{3}_{-6}\otimes E^{1}_{-1}\to E^{4}_{-7}, \\
E^{2}_{-3}\otimes E^{2}_{-4}\to E^{4}_{-7},\qquad 
E^{2}_{-4}\otimes E^{2}_{-3}\to E^{4}_{-7}.
\end{gather*}
By Theorem \ref{yythm1.8}, if $A$ is an AS regular algebra, then $E$ is a
Frobenius algebra. The multiplication $m_2$ of a Frobenius algebra $E$ of 
type \type{12221} can be described as follows. 

Let $\delta$ be a basis element of $E^{4}_{-7}$. Pick a basis element
$\beta_1\in E^{2}_{-3}$. Since $E$ is a Frobenius algebra,
\[
m_2: E^{2}_{-3}\otimes E^{2}_{-4}\to E^{4}_{-7} \quad \text{and} \quad 
m_2: E^{2}_{-4}\otimes E^{2}_{-3}\to E^{4}_{-7}
\]
are both nonzero. So we can pick a basis element $\beta_2\in E^{2}_{-4}$
such that $\beta_1 \beta_2=\delta$ and $\beta_2 \beta_1=t\; \delta$ for some 
$0\neq t \in k$. Pick a basis $\{\alpha_1, \alpha_2\}$ for $E^{1}_{-1}$. 
Since $E$ is a Frobenius  algebra, 
\[
m_2: E^{1}_{-1}\otimes E^{3}_{-6}\to E^{4}_{-7} \quad \text{and} \quad 
m_2: E^{3}_{-6}\otimes E^{1}_{-1}\to E^{4}_{-7}
\]
are perfect pairings. Hence we may choose a basis $\{\gamma_1, \gamma_2\}$ of
$E^{3}_{-6}$ such that $\alpha_i \gamma_j=\delta_{ij}\delta$.
Let $\gamma_i \alpha_j=r_{ij}\delta$ for some $r_{ij}\in k$. Then the matrix
$\Lambda:=(r_{ij})_{2\times 2}\in M_2(k)$ is non-singular. 

We call $(\Lambda, t)$ the \emph{Frobenius data} of $E$ or of $A$. 

Possible nonzero components of $m_3$ on $E^{\otimes 3}$ are
\begin{gather*}
(E^{1}_{-1})^{\otimes 3}\to E^{2}_{-3}, \\
(E^{1}_{-1})^{\otimes 2}\otimes E^{2}_{-4}\to E^{3}_{-6},\quad
E^{1}_{-1}\otimes E^{2}_{-4}\otimes E^{1}_{-1}\to E^{3}_{-6},\quad
E^{2}_{-4}\otimes 
(E^{1}_{-1})^{\otimes 2}\to E^{3}_{-6}, \\
E^{1}_{-1}\otimes (E^{2}_{-3})^{\otimes 2}\to E^{4}_{-7},\quad
E^{2}_{-3}\otimes E^{1}_{-1}\otimes E^{2}_{-3}\to E^{4}_{-7}, \quad 
(E^{2}_{-3})^{\otimes 2}\otimes E^{1}_{-1}\to E^{4}_{-7}.
\end{gather*}
We have, for
$1 \leq i, j, k \leq 2$,
\begin{gather*}
m_3(\alpha_i,\alpha_j,\alpha_k)=a_{ijk} \; \beta_1; \\
m_3(\alpha_i,\alpha_j,\beta_2)=b_{13ij} \; \gamma_1+b_{23ij} \; \gamma_2,
\\
m_3(\alpha_i,\beta_2,\alpha_j)=b_{12ij} \; \gamma_1+b_{22ij} \; \gamma_2,
\\
m_3(\beta_2,\alpha_i,\alpha_j)=b_{11ij} \; \gamma_1+b_{21ij} \; \gamma_2;
\\
m_3(\alpha_i,\beta_1,\beta_1)=c_{1i} \; \delta, 
\\
m_3(\beta_1,\alpha_i,\beta_1)=c_{2i} \; \delta, 
\\
m_3(\beta_1,\beta_1,\alpha_i)=c_{3i} \; \delta,
\\
\text{all other applications of $m_3$ are zero.}
\end{gather*}
Here $a_{ijk}$, $b_{ipjk}$ and $c_{pi}$ are scalars in the field $k$. 

Since $m_1=0$, the Stasheff identity \SI{4} (see Definition \ref{yydef2.1})
becomes 
\[
m_3(m_2\otimes id^{\otimes 2}-id\otimes m_2\otimes id+id^{\otimes 2}
\otimes m_2)-m_2(m_3\otimes id+id\otimes m_3)=0.
\]
Applying \SI{4} to the elements $(\alpha_i,\alpha_j,\alpha_k,\beta_2)$,
$(\alpha_i,\alpha_j,\beta_2,\alpha_k)$, $(\alpha_i,\beta_2,\alpha_j,\alpha_k)$,
and $(\beta_2,\alpha_i,\alpha_j,\alpha_k)$, respectively, we obtain
\begin{equation}
\begin{gathered}
a_{ijk} = b_{i3jk}, \quad 
b_{i2jk}= \sum_{s=1}^2 r_{sk}\; b_{s3ij}, \\
b_{i1jk}= \sum_{s=1}^2 r_{sk}\; b_{s2ij}, \quad 
-t\; a_{ijk} = \sum_{s=1}^2 r_{sk}\; b_{s1ij}.
\end{gathered}
\tag*{\SI{4a}}
\end{equation}
Note that Koszul sign rule applies when two symbols are commuted.
As a consequence of \SI{4a}, we have
\begin{equation}
-t\; a_{ijk}=
\sum_{s,t,u=1}^2r_{sk}\; r_{tj}\; r_{ui}\; a_{uts}.
\tag*{\SI{4b}}
\end{equation}

Possible nonzero applications of $m_4$ on $E^{\otimes 4}$ are
\begin{gather*}
(E^{1}_{-1})^{\otimes 4}\to E^{2}_{-4},
\\
(E^{1}_{-1})^{\otimes 3}\otimes E^{2}_{-3}\to E^{3}_{-6},\qquad
(E^{1}_{-1})^{\otimes 2}\otimes E^{2}_{-3}\otimes E^{1}_{-1}\to E^{3}_{-6}, 
\\
E^{1}_{-1}\otimes E^{2}_{-3}\otimes (E^{1}_{-1})^{\otimes 2}\to E^{3}_{-6},
\qquad
E^{2}_{-3}\otimes (E^{1}_{-1})^{\otimes 3}\to E^{3}_{-6}.
\end{gather*}
We write down the coefficients of these maps.  For
$1 \leq i, j, k, h \leq 2$,
\begin{gather*}
m_4(\alpha_i,\alpha_j,\alpha_k,\alpha_h)=y_{ijkh} \;  \beta_2,
\\
m_4(\alpha_i,\alpha_j,\alpha_k,\beta_1)=
x_{14ijk}\;\gamma_1+x_{24ijk}\;\gamma_2,
\\
m_4(\alpha_i,\alpha_j,\beta_1,\alpha_k)=
x_{13ijk}\;\gamma_1+x_{23ijk}\;\gamma_2,
\\
m_4(\alpha_i,\beta_1,\alpha_j,\alpha_k)=
x_{12ijk}\;\gamma_1+x_{22ijk}\;\gamma_2,
\\
m_4(\beta_1,\alpha_i,\alpha_j,\alpha_k)=
x_{11ijk}\;\gamma_1+x_{21ijk}\;\gamma_2,
\end{gather*}
$$\text{all other applications of $m_4$ are zero.}$$
Here $y_{ijkh}$ and $x_{hpijk}$ are scalars in the field $k$. 

The Stasheff identity \SI{5} becomes
\begin{multline*}
m_4(m_2\otimes id^{\otimes 3}-id\otimes m_2\otimes id^{\otimes 2}
+id^{\otimes 2}\otimes m_2\otimes id-id^{\otimes 3}\otimes m_2)+ \\
+m_3(m_3\otimes id^{\otimes 2}+id\otimes m_3\otimes id+id^{\otimes 2}\otimes
m_3)+m_2(m_4\otimes id-id\otimes m_4)=0.
\end{multline*}
As with \SI{4}, after applied to elements 
$(\alpha_i,\alpha_j,\alpha_k,\alpha_h,\beta_1)$, 
$(\alpha_i,\alpha_j,\alpha_k,\beta_1,\alpha_h)$, 
$(\alpha_i,\alpha_j,\beta_1,\alpha_k,\alpha_h)$, 
$(\alpha_i,\beta_1,\alpha_j,\alpha_k,\alpha_h)$ and  
$(\beta_1,\alpha_i,\alpha_j,\alpha_k,\alpha_h)$,
respectively, \SI{5} gives the following equations:
\begin{equation}
\begin{gathered}
a_{ijk} \; c_{2h}-a_{jkh} \; c_{1i}+t \; y_{ijkh}
-x_{i4jkh} =0, \\
a_{ijk} \; c_{3h}+r_{1h}\; x_{14ijk} +r_{2h}\;x_{24ijk} 
-x_{i3jkh} =0, \\
r_{1h}\; x_{13ijk} +r_{2h}\; x_{23ijk} 
-x_{i2jkh} =0, \\
c_{1i}\; a_{jkh} -r_{1h}\; x_{12ijk} -r_{2h}\; x_{22ijk} 
+x_{i1jkh} =0, \\
a_{jkh} \; c_{2i}-a_{ijk} \; c_{3h}
-r_{1h}\; x_{11ijk} -r_{2h}\; x_{21ijk} + y_{ijkh}=0.
\end{gathered}
\tag*{\SI{5a}}
\end{equation}

The Stasheff identity \SI{6} becomes
\begin{multline*}
\quad \quad m_4(-m_3\otimes id^{\otimes 3}-id\otimes m_3\otimes id^{\otimes 2}
-id^{\otimes 2}\otimes m_3\otimes id-id^{\otimes 3}\otimes m_3) \\
+m_3(m_4\otimes id^{\otimes 2}-id\otimes m_4\otimes id
+id^{\otimes 2}\otimes m_4)=0. \quad \quad 
\end{multline*}
Applying \SI{6} to $(\alpha_i,\alpha_j,\alpha_k,\alpha_h,\alpha_m,\alpha_n)$,
we obtain the equation
\begin{multline*}
\quad \quad -a_{ijk}x_{s1hmn}+a_{jkh}x_{s2imn}-a_{khm}x_{s3ijn}+a_{hmn}x_{s4ijk} \\
+b_{s1mn}y_{ijkh}-b_{s2in}y_{jkhm}+b_{s3ij}y_{khmn}=0. \quad \quad 
\tag*{\SI{6a}}
\end{multline*}
We now have all of the equations we need, by the following lemma.

\begin{lemma}
\label{yylem4.1} 
If $E=E(m_{2}, m_{3}, m_{4})$ is as in Lemma~\ref{yylem3.2}(a), then
the identity \SI{n} holds for every $n\geq 7$.
\end{lemma}

\begin{proof} 
In this situation, the Stasheff identity \SI{7} becomes
\[
m_4(\sum_{i=0}^{3}(-1)^{i}id^{\otimes i}\otimes m_4\otimes id^{\otimes 3-i})=0,
\]
which is zero when applied to basis elements. For $n\geq 8$, 
\SI{n} must involve $m_i$ for some $i\geq 5$, and hence it is zero.
\end{proof}

Given an AS regular algebra $A=k\langle z_1,z_2\rangle/(R)$ of type \type{12221}, 
we know from Proposition \ref{yyprop1.4} that $A$ has two relations $r_3$ 
and $r_4$ of degree 3 and 4, respectively, which we write as
\begin{equation}
\label{E4.1.1}
r_3=\sum a_{ijk}\; z_iz_jz_k
\end{equation}
and
\begin{equation}
\label{E4.1.2}
r_4=\sum y_{ijkh}\; z_iz_jz_kz_h.
\end{equation}
By Theorem \ref{yythm2.2}, for the Ext-algebra $E$ of $A$, we have
\begin{equation}\label{E4.1.3}
\begin{gathered}
m_3(\alpha_i,\alpha_j,\alpha_k)=a_{ijk}\; \beta_1 \\
m_4(\alpha_i,\alpha_j,\alpha_k,\alpha_h)=y_{ijkh}\; \beta_2.
\end{gathered}
\end{equation}
Conversely, if we know for the Ext-algebra $E$ that \eqref{E4.1.3}
holds, then the relations of $A$ are given by 
\eqref{E4.1.1} and \eqref{E4.1.2}. The main idea of the next section 
is to classify all possible higher multiplications $m_3$ and $m_4$ on
$E$, thus all possible coefficients $a_{ijk}$ and $y_{ijkh}$.  Then
we define the relations $r_3$ and $r_4$ using \eqref{E4.1.1} 
and \eqref{E4.1.2}, and we investigate when the resulting algebra 
$k\langle z_1,z_2\rangle/(r_3,r_4)$ is AS regular.

By Lemma \ref{yylem1.3} if $A$ is an AS regular algebra of type 
\type{12221}, then $A$ is a domain.  So the relations satisfy the following
two conditions.
\begin{equation}
\label{E4.R3}
\text{$r_3$ is neither zero nor a product of lower degree polynomials,}
\tag{R3}
\end{equation}
and
\begin{equation}
\label{E4.R4}
\text{$r_4$ is neither zero nor a product of lower degree polynomials.}
\tag{R4}
\end{equation}

By Proposition \ref{yyprop1.4}(a), the first ten terms of the Hilbert
series of $A$ are
\begin{equation}
\label{E4.HS}
H_A(t)
=1+2t+4t^2+7t^3+11t^4+16t^5
+23 t^6+31t^7+41t^8+53t^9+67t^{10}+\cdots.
\tag{HS}
\end{equation}

\section{Generic algebras of type \type{12221}}
\label{yysect5}

In this section we classify generic AS regular algebras of type
\type{12221}. The idea is simple. We first classify all generic
multiplications $m_2$ on $E$, then all generic maps $m_3$, and then
all possible maps $m_4$. After we know $m_3$ and $m_4$, we use Theorem
\ref{yythm2.2} to recover the relations $r_3$ and $r_4$ of $A$.

In this section let $k$ be an algebraically closed field,
and let $A$ be an AS regular algebra of type \type{12221}.

We first work on $m_2$ of $E$.  If $B$ is an invertible $2 \times 2$
matrix, then replacing $(\alpha_1,\alpha_2)$ by $(\alpha_1,
\alpha_2)B^{-1}$ and $(\gamma_1,\gamma_2)$ by $(\gamma_1,
\gamma_2)B^{T}$ changes the Frobenius data of $E$ from $(\Lambda,t)$
to $(B\Lambda B^{-1},t)$.  So choosing $B$ properly, we may assume
that the matrix $\Lambda$ is either
\[
\begin{pmatrix} g_1 &0\\ 0&g_2\end{pmatrix} \qquad \text{or}
\qquad \begin{pmatrix} g_1 &1\\ 0&g_1\end{pmatrix}
\]
for some nonzero $g_i\in k$. 

We now give a generic condition for $m_2$:
\begin{equation}
\text{Let $g_1$ and $g_2$ be the eigenvalues of $\Lambda$.  Then
$(g_1g_2^{-1})^i\neq 1$ for $1\leq i\leq 4$.}
\tag{GM2}
\end{equation}

\begin{remark}
\label{yyrem5.1} Let $E$ be the Ext-algebra of an AS regular algebra of
type \type{12221}.  Let $(\Lambda,t)$ be the Frobenius data of $E$,
and let $g_i$ be the eigenvalues of $\Lambda$.
By some preliminary results in \cite{LWZ}, the following are equivalent:
\begin{enumerate}
\item $g_1$ is a root of 1.
\item $g_2$ is a root of 1.
\item $t$ is a root of 1.
\item $g_1g_2^{-1}$ is a root of 1.
\end{enumerate}
\end{remark}

Ideally if $A$ is ``generic'', then $g_i$ should not be a root of
unity.  By Remark \ref{yyrem5.1}, this implies that $g_1g_2^{-1}$ is
not a root of 1, so (GM2) holds.  Hence (GM2) can be viewed as a
sort of generic condition.

Suppose now that (GM2) holds.  Then $g_1\neq g_2$, so we may assume that 
$\Lambda=\begin{pmatrix} g_1 &0\\ 0&g_2\end{pmatrix}$.  The
multiplication map $m_2$ is now described by the Frobenius data. 

Next we work with $m_3$ by considering \SI{4}. By \SI{4b}, we have
\[
-t a_{ijk}=g_ig_jg_k a_{ijk} \quad \text{for all}\quad i,j,k=1,2,
\]
or equivalently, 
\begin{gather*}
(t+\; g_1^3\; )\; a_{111}\; =0,
\\
(t+g_1^2g_2)\; a_{112}=0, \quad 
(t+g_1^2g_2)\; a_{121}=0, \quad 
(t+g_1^2g_2)\; a_{211}=0,
\\
(t+g_1g_2^2)\; a_{122}=0, \quad 
(t+g_1g_2^2)\; a_{212}=0, \quad 
(t+g_1g_2^2)\; a_{221}=0,
\\
(t+\; g_2^3\; )\; a_{222}\; =0.
\end{gather*}
If all $t+g_1^s g_2^{3-s}\neq 0$ for $s=0,1,2,3$, then the only
solution is $a_{ijk}=0$ for all $i,j,k$.  In this case, the corresponding $r_3$
has zero coefficients, a contradiction. It follows from (GM2) that at
most one of the $t+g_1^sg_2^{3-s}$ terms can be zero. If $t+g_1^3=0$,
then $a_{ijk}=0$ except for $a_{111}$; thus $r_3=a_{111}z_1^3$, which
does not satisfy \eqref{E4.R3}.  So this is impossible.  The same
holds for the case $t+g_2^3=0$.  Therefore we have either
\[
t+g_1^2g_2=0 \quad\text{or}\quad t+g_1g_2^2=0.
\]
By switching $\alpha_1\leftrightarrow \alpha_2$ and 
$\gamma_1\leftrightarrow \gamma_2$, these two cases are equivalent. 
Therefore, without loss of generality,
we may assume that
\[
t+g_1g_2^2=0.  
\]
Under this hypothesis, the coefficients $a_{ijk}$ are zero except
possibly for $a_{122}$, $a_{212}$, and $a_{221}$.  Hence 
\[
r_3=a_{122}z_1z_2^2+a_{212}z_2z_1z_2+a_{221}z_2^2z_1.
\]
Since $A$ is a domain, it follows from \eqref{E4.R3} that $a_{122}\neq 0$
and $a_{221}\neq 0$. Without changing the structure of $E$ we may assume 
that 
\begin{equation}\label{E5.A}
\begin{gathered}
a_{122}=1,\;  a_{212}=v,\;  a_{221}=w\neq 0, \\
a_{ijk}=0\quad \text{for all other $(ijk)$}.
\end{gathered}
\end{equation}
Using \SI{4a}, we solve for $b_{iqjk}$; the solution is
\begin{equation}
\label{E5.B}
\begin{gathered}
b_{iqjk}=0 \quad \text{except for}  \\
\begin{matrix}
b_{1322} = 1, & b_{2312} = v, & b_{2321} = w, \\
b_{2221} = g_{1}, & b_{1222} = g_{2}v, & b_{2212} = g_{2}w,\\ 
b_{1122} = g_{2}^2w, & b_{2112} = g_{2}g_{1}, & b_{2121} = g_{1}g_{2}v.
\end{matrix}
\end{gathered}
\end{equation}

Our generic condition for $m_3$ is 
\begin{equation}
a_{122}+a_{212}+a_{221}\neq 0 \quad \text{or} \quad 1+v+w\neq 0.
\tag{GM3}
\end{equation}
This condition is not very essential, but it guarantees that the matrix
$\Lambda$ is diagonal.  An AS regular algebra of type \type{12221} is 
called \emph{$(m_2,m_3)$-generic} if its $A_\infty$-Ext-algebra satisfies
both (GM2) and (GM3). 

\begin{lemma}
\label{yylem5.2} Assume that $E$ satisfies \textup{\eqref{E5.A}} and 
\textup{(GM3)}. Then 
$\Lambda=\begin{pmatrix} g_1 &0\\ 0&g_2\end{pmatrix}$ and $t=-g_1g_2^2$.
\end{lemma}

\begin{proof} Let $g_1$ and $g_2$ be the eigenvalues of $\Lambda=
\begin{pmatrix} r_{11} &r_{12}\\ r_{21}&r_{22}\end{pmatrix}$.  A direct
(but tedious) computation from \SI{4b} together with \eqref{E5.A} and (GM3)
shows that $r_{12}=0=r_{21}$ and $t=-r_{11}r_{22}^2$.
\end{proof}

Without (GM3), we cannot conclude that $r_{12}=0$.  Lemma \ref{yylem5.2}
means that under the hypotheses of \eqref{E5.A} and (GM3),
the matrix $\Lambda$ is always diagonal, so when we make changes to
the relation $r_4$, the structure of $m_2$ and $m_3$ will not change. 
By replacing $r_4$ by the equivalent relation
\[
r_4-y_{1122}z_1r_3-y_{2122}z_2r_3-y_{1221}r_3z_1-y_{1222}r_3z_2,
\]
we may assume that
\begin{equation}
\label{E5.2.1}
y_{1122}=y_{2122}=y_{1221}=y_{1222}=0.
\end{equation}
This is the advantage of (GM3).  Using \SI{5a} recursively, we obtain
(when $\Lambda$ is diagonal) that 
\begin{equation}
\label{E5.SI(5-6)}
(1-tg_{i}g_{j}g_{k}g_{h})y_{ijkh}=a_{ijk}U_{ijkh}+a_{jkh}V_{jkhi}
\tag*{\SI{5b}}
\end{equation}
for all $i,j,k$. Here
\begin{gather*}
U_{ijkh}=-g_{h}c_{1h}+g_{i}g_{j}g_{k}g_{h}c_{2h}+c_{3h}, \\
V_{jkhi}=-g_{i}g_{j}g_{k}g_{h}c_{1i}-c_{2i}+g_{j}g_{k}g_{h}c_{3i}.
\end{gather*}
Using \eqref{E5.A} and \eqref{E5.2.1}, the equations \ref{E5.SI(5-6)}
become
\begin{align*}
(1-tg_{1}^4)y_{1111} &= 0, &  (1-tg_{2}^4)y_{2222} &= 0, \\
(1-tg_{2}g_{1}^3)y_{1112} &= 0, &  (1-tg_{2}^2g_{1}^2)y_{1212} &= vV_{2121}, \\
(1-tg_{2}g_{1}^3)y_{1121} &= 0, &  (1-tg_{2}^2g_{1}^2)y_{2211} &= wU_{2211}, \\
(1-tg_{2}g_{1}^3)y_{1211} &= 0, &  (1-tg_{2}^2g_{1}^2)y_{2121} &= vU_{2121}, \\
(1-tg_{2}g_{1}^3)y_{2111} &= 0, &  (1-tg_{2}^2g_{1}^2)y_{2112} &= 0,         \\
(1-tg_{2}^3g_{1})y_{2212} &= wU_{2212}+vV_{2122}, & 
(1-tg_{2}^3g_{1})y_{2221} &= wV_{2212}, \\
V_{1221} &= 0,  &  vU_{2122}+V_{1222} &= 0, \\
U_{1222} &= 0,  &  U_{1221}+wV_{2211} &= 0.
\end{align*}
Since $U_{ijkh}=U_{i'j'k'h}$ and $V_{ijkh}=V_{i'j'k'h}$ when
the set $\{i,j,k\}$ is equal to the set $\{i',j',k'\}$,
the last four equations imply
\begin{gather*}
U_{1221}=U_{2211}=U_{2121}=0, \\
U_{1222}=U_{2122}=U_{2212}=0, \\
V_{1221}=V_{2211}=V_{2121}=0, \\
V_{1222}=V_{2122}=V_{2212}=0. 
\end{gather*}
Hence \ref{E5.SI(5-6)} becomes
\begin{equation}\label{E5.SI(5-7)}
(1-tg_ig_jg_kg_h)y_{ijkh}=0
\tag*{\SI{5c}}
\end{equation}
for all $i,j,k,h$.  If all of the terms $1-tg_1^sg_2^{4-s}$ are
nonzero, then all $y_{ijkh}=0$, and this contradicts \eqref{E4.R4}.
The generic condition (GM2) says that $(g_1g_2^{-1})^i\neq 1$ for
$1\leq i\leq 4$; hence at most one of $1-tg_1^sg_2^{4-s}$ is zero.  Of
course, at least one must be zero; otherwise $r_4=0$, which contradicts
\eqref{E4.R4}.  We consider the possibilities in the following.

\begin{proposition}\label{prop5.new}
Suppose that $A$ is an AS regular algebra satisfying the generic
conditions (GM2) and (GM3).  
\begin{description}
\item[Case 1] $1-tg_1^4=0$.  There are no AS regular algebras
satisfying this.
\item[Case 2] $1-tg_2^4=0$.  There are no AS regular algebras
satisfying this.
\item[Case 3] $1-tg_1g_2^3=0$.  There are no AS regular algebras
satisfying this.
\item[Case 4] $1-tg_1^2g_2^2=0$.  There are no AS regular algebras
satisfying this.
\item[Case 5] $1-tg_1^3g_2=0$.  This leads to several possible
relations $r_{3}$ and $r_{4}$, all of which give AS regular algebras.
\end{description}
\end{proposition}

The proofs for cases 4 and 5 are long.  We deal with cases 1--3 here,
and defer cases 4 and 5 to subsequent sections.

\begin{proof}
Cases 1 and 2 contradict \eqref{E4.R4}: in either of these cases, the
relation $r_{4}$ is a product of lower degree polynomials.  The same
holds for Case 3, because of assumption \eqref{E5.2.1}.
\end{proof}

We will analyze Case 4 in Section \ref{yysect6}, where we will show
that those algebras are not AS regular: see
Lemmas~\ref{yylem6.2}, \ref{yylem6.3}, and \ref{yylem6.4}.  We will
analyze Case 5 in Section~\ref{yysect-reg}.

\section{Case 4: Non-regular algebras}
\label{yysect6}

We continue our discussion of AS regular algebras of type \type{12221}
which satisfy the two generic conditions (GM2) and (GM3).  In this
section we assume that we are in Case 4 of
Proposition~\ref{prop5.new}.  We show that the resulting algebras
cannot in fact be AS regular.

We have an algebra $A$ of the form 
\[
A = k\langle z_1,z_2\rangle/(r_3,r_4),
\]
where $k$ is a field, and the relations are
\begin{gather*}
r_3= z_1z_2^2+vz_2z_1z_2+wz_2^2z_1,\\ 
r_4=\sum y_{ijkh}\; z_iz_jz_kz_h.
\end{gather*}
According to equation~\ref{E5.SI(5-7)} and the discussion following
it, Case 4 says that the possible nonzero $y$'s are
\[
y_{1212},y_{2112},y_{2121},y_{2211}.
\]
By \eqref{E4.R4}, $y_{1212}\neq 0$.  So we may assume that
\[
y_{1212}=1,y_{2112}=p,y_{2121}=q,y_{2211}=r, \quad \text{and}\quad
y_{ijkh}=0
\]
for all other $(i,j,k,h)$.  That is, we may assume that the relations
for $A$ have the form
\begin{gather*}
r_3= z_1z_2^2+vz_2z_1z_2+wz_2^2z_1, \\
r_4=z_1z_2z_1z_2+pz_2z_1^2z_2+qz_2z_1z_2z_1+rz_2^2z_1^2.
\end{gather*}
In this case we may assume that $c_{i2}=0$ for all $i$, by the next lemma.

\begin{lemma}
\label{yylem6.1}
Suppose $A=k\langle z_1,z_2\rangle/(r_3,r_4)$, where
$r_3= z_1z_2^2+vz_2z_1z_2+wz_2^2z_1$ and $r_4=z_1z_2z_1z_2+pz_2z_1^2z_2+
qz_2z_1z_2z_1+rz_2^2z_1^2$. If $A$ is an AS regular algebra of dimension
4, then the coefficient $c_{i2}$ in $E$ is zero for all $i$.
\end{lemma}

\begin{proof} 
We define an Adams grading on $A$ by
\[
\adeg(z_1)=(1,1,0)\in \Z^3, \quad
\adeg(z_2)=(1,0,1)\in \Z^3.
\]
The first component gives the usual connected grading of $A$.  Note
that $\adeg(r_3)=(3,1,2)$ and $\adeg(r_4)=(4,2,2)$.  The Ext-algebra
$E$ will also have an Adams $\Z^3$-grading, with
\[
\adeg(\alpha_1)=(-1,-1,0) \quad \text{and} \quad
\adeg(\alpha_2)=(-1,0,-1).
\]
From this, together with \eqref{E4.1.3} and the
definition of $\delta$, it is easy to compute the following Adams
degrees:
\begin{gather*}
\adeg(\beta_1)=(-3,-1,-2), \quad 
\adeg(\beta_2)=(-4,-2,-2), \quad \text{and} \\
\adeg(\delta)=(-7,-3,-4).
\end{gather*}
Also
\[
\adeg(\alpha_2\otimes\beta_1\otimes\beta_1)=(-7,-2,-5).
\]
Hence $c_{i2}=0$ by its definition in terms of $m_3$.
\end{proof}

Next we perform some long but elementary computations.  There are
similar computations for Case 5, and we give a few more details for
those; see the discussion after Lemma~\ref{yylem5.3}.

We start by using \SI{5} to find formulas for $x_{isjkh}$ (which we
omit), and then we input those into \SI{6}.  This produces $2^7$
equations, which can be generated by Maple 8.  After a few steps of
simplification, we are able to list all possible solutions.  The
details are omitted since the computations are straightforward; we
only give the solutions here.  After deleting some useless solutions
(they give algebras that are clearly not domains), we have three
non-trivial solutions, as listed below.

Set $f=c_{21}$ and $h=g_{1}g_{2}$.  Note that $h\neq 0$.  To make
sure that the algebra is a domain, we also need that $w\neq 0$.

\begin{solution}\label{sol2.1}
\begin{gather*}
v = h^2 f-h,\quad w = -h^3 f,
\\
p = h^2 f,\quad q = h^3 f,\quad  r = h^5 f^2.
\end{gather*}
\end{solution}

\begin{solution}\label{sol2.2}
\begin{gather*}
v = 0,\quad w = -h^2,
\\
p = h,\quad q = h^2, \quad r = h^4 f.
\end{gather*}
\end{solution}

\begin{solution}\label{sol2.3}
This is the most complicated one. Since $h\neq 0$, $h$ is invertible.
We find that $w=-hp$. Since $w\neq 0$, we have $p\neq 0$
and hence $p$ is also invertible. When $h\neq p$, the solution is
\begin{gather*}
v = h^{-1}p^{-1}(h^3-p^3),\quad w = -h p,
\\
p=p,\quad 
q = h^{-2}p^{-1}(h^5+h^4 p+h^3 p^2-h p^4-p^5)\quad 
r = -p (h^2+h p+p^2).
\end{gather*}
In this case $f=-h^{-3}(h^2+hp+p^2)$. When $h=p$, this reduces to 
Solution~\ref{sol2.2}.
\end{solution}

We then use Bill Schelter's program ``Affine'' to check the Hilbert series
of the algebras in the above three cases.

The algebra corresponding to Solution~\ref{sol2.1} is $X(p,h)=k\langle z_1,z_2
\rangle/(r_3,r_4)$, where the two relations are  
\begin{gather*}
r_3= z_1z_2^2+(p-h)z_2z_1z_2-hpz_2^2z_1, \\ 
r_4=z_1z_2z_1z_2+pz_2z_1^2z_2+hpz_2z_1z_2z_1+hp^2z_2^2z_1^2.
\end{gather*}
If either $p$ or $h$ is zero, then $X(p,h)$ is not a domain.

\begin{lemma}
\label{yylem6.2} The algebra $X(p,h)$ is not a noetherian 
AS regular algebra of dimension 4.
\end{lemma}

\begin{proof} Let $X=X(p,h)$. We prove the assertion by showing 
that the Hilbert series of $X$ is wrong.

The Affine program gives the first few terms of the 
Hilbert series of $X$:
\[
H_X(t)=1+2t+4t^2+7t^3+11t^4+17t^5+\cdots,
\]
which is different from \eqref{E4.HS}: in \eqref{E4.HS} the 
coefficient of $t^5$ is 16. 
\end{proof}

There are other ways to prove this lemma without 
using Affine. For example, the $k$-vector space dimension 
of $X_5$ can be computed directly by using the ideas in Bergman's 
diamond lemma \cite{Be}. It can also be proved that the Hilbert 
series of $X$ (for any $p$ and $h$) is equal to the Hilbert 
series of $B:=k\langle z_1,z_2\rangle/(z_1z_2^2,z_1z_2z_1z_2)$
(which is $X(0,0)$ by definition). Now $B$ is a monomial algebra,
and the start of $H_{B}(t)$ is $1+2t+4t^2+7t^3+11t^4+17t^5$
by a direct computation. 

The algebra corresponding to Solution~\ref{sol2.2} is $Y(h,f)=k\langle z_1,z_2
\rangle/(r_3,r_4)$, where 
\begin{gather*}
r_3= z_1z_2^2-h^2z_2^2z_1, \\
r_4=z_1z_2z_1z_2+hz_2z_1^2z_2+h^2z_2z_1z_2z_1+h^4fz_2^2z_1^2.
\end{gather*}

\begin{lemma}
\label{yylem6.3} The algebra $Y(h,f)$ is not a noetherian 
AS regular algebra of dimension 4.
\end{lemma}

\begin{proof} Let $Y=Y(h,f)$. Suppose to the contrary that $Y$ is a noetherian 
AS regular algebra of dimension 4.  By the relation $r_3$, $z_2^2$ is 
normal element.  By Lemma~\ref{yylem1.3}, $Y$ is a domain, so $z_2^2$ 
is regular.  Let $B=Y/(z_2^2)$; then $H_B(t)=
(1-t^2)H_Y(t)$.  It follows from Proposition \ref{yyprop1.4}(a) that 
\begin{equation}
\label{E6.3.1}
H_B(t)=\frac{1}{(1-t)^2(1-t^3)}=1+2t+3t^2+5t^3+7t^4+9t^5
+\cdots.
\end{equation}
Clearly $B$ is isomorphic to $k\langle z_1,z_2
\rangle/(z_2^2,z_1z_2z_1z_2+hz_2z_1^2z_2+h^2z_2z_1z_2z_1)$,
which is independent of the parameter $f$.
The Affine program gives the following Hilbert series for $B$
(for any $h$):
\[
H_B(t)=1+2t+3t^2+5t^3+7t^4+10t^5+
\cdots.
\]
This is different from \eqref{E6.3.1}, yielding a contradiction. 
\end{proof}

Another way of proving the above lemma is to show that the 
Hilbert series of $B$ is equal to the Hilbert series of 
$B':=k\langle z_1,z_2\rangle/(z_2^2,z_1z_2z_1z_2)$. 
Since $B'$ is a monomial algebra, the computation of the Hilbert 
series of $B'$ is straightforward. 

We can also compute the Hilbert series of $Y$ directly 
without introducing $B$. The Affine program gives
\[
H_Y(t)=1+2t+4t^2+7t^3+11t^4+16t^5+23t^6+32t^7+\cdots.
\]
The coefficient of the $t^7$ term is different from that in \eqref{E4.HS};
hence the algebra $Y$ is not AS regular.  The direct computation 
of the Hilbert series of $Y$ without the Affine program is
tedious. 

The algebra corresponding to Solution~\ref{sol2.3} is $Z(p,h)=k\langle z_1,z_2
\rangle/(r_3,r_4)$, where 
\begin{align*}
r_3 &= z_1z_2^2+h^{-1}p^{-1}(h^3-p^3)z_2z_1z_2-hpz_2^2z_1, \\
r_4 &= z_1z_2z_1z_2+pz_2z_1^2z_2-p(h^2+hp+p^2)z_2^2z_1^2 \\
& \quad + h^{-2}p^{-1}(h^5+h^4p+h^3p^2-hp^4-p^5)z_2z_1z_2z_1.
\end{align*}
The relations here are more complicated than last two cases. 
We handle this algebra in the same fashion as the last case.
For example, for any nonzero $p$ and $h$ except for the case
$p+h=0$, the Affine program gives the following Hilbert series:
\[
H_{Z(p,h)}(t)=1+2t+4t^2+7t^3+11t^4+16t^5+23t^6+32t^7+\cdots
\]
which is different from \eqref{E4.HS}.  We would like to remark
that the higher-degree terms of $H_{Z(p,h)}(t)$ are different 
for different values of $(p,h)$, but for our purposes, the 
coefficient of $t^7$ provides enough information.  When $p+h=0$, 
the Affine program gives the following Hilbert series:
\[
H_{Z(p,-p)}(t)=1+2t+4t^2+7t^3+11t^4+17t^5+26t^6+39t^7+\cdots
\]
which is also different from \eqref{E4.HS}.
Hence the following holds.

\begin{lemma}
\label{yylem6.4} The algebra $Z(p,h)$ is not a noetherian 
AS regular algebra of dimension 4.
\end{lemma}

\section{Case 5: Regular algebras}
\label{yysect-reg}

We continue our discussion of AS regular algebras of type \type{12221}
which satisfy the two generic conditions (GM2) and (GM3).  In this
section we assume that we are in Case 5 of
Proposition~\ref{prop5.new}.  We prove that these algebras are indeed
AS regular, and we describe some of their properties.

We have an algebra $A$ of the form 
\[
A = k\langle z_1,z_2\rangle/(r_3,r_4),
\]
where $k$ is an algebraically closed field, and the relations are
\begin{gather*}
r_3= z_1z_2^2+vz_2z_1z_2+wz_2^2z_1, \\
r_4=\sum y_{ijkh}\; z_iz_jz_kz_h.
\end{gather*}
According to equation~\ref{E5.SI(5-7)} and the discussion following
it, Case 5 says that the possible nonzero $y$'s are
\[
y_{1112},y_{1121},y_{1211},y_{2111}.
\]
By \eqref{E4.R4}, $y_{1112}\neq 0$ and $y_{2111}\neq 0$, and we may
assume that
\[
y_{1112}=1, \ y_{1121}=p, \ y_{1211}=q, \ y_{2111}=r\neq 0,
\quad\text{and}\quad y_{ijkh}=0
\]
for all other $(i,j,k,h)$.  That is, we may assume that the relations
for $A$ have the form
\begin{gather*}
r_3= z_1z_2^2+vz_2z_1z_2+wz_2^2z_1, \\
r_4=z_1^3z_2+pz_1^2z_2z_1+ qz_1z_2z_1^2+rz_2z_1^3.
\end{gather*}
In this case we can show that $c_{ih}=0$ for all $i,h$.

\begin{lemma}
\label{yylem5.3} Suppose $A=k\langle z_1,z_2\rangle/(r_3,r_4)$ where
$r_3= z_1z_2^2+vz_2z_1z_2+wz_2^2z_1$ and $r_4=z_1^3z_2+pz_1^2z_2z_1+
qz_1z_2z_1^2+rz_2z_1^3$. If $A$ is an AS regular algebra of dimension
4, then the coefficient $c_{ih}$ in $E$ is zero.
\end{lemma}

\begin{proof} 
Just as in the proof of Lemma \ref{yylem6.1}, we define an Adams
grading on $A$ by
\[
\adeg(z_1)=(1,1,0)\in \Z^3, \quad
\adeg(z_2)=(1,0,1)\in \Z^3.
\]
It is easy to compute the following Adams degrees:
\[
\adeg(\beta_1)=(-3,-1,-2),\quad \adeg(\delta)=(-7,-4,-3),
\]
and
\[
\adeg(\alpha_1\otimes\beta_1\otimes\beta_1)=(-7,-3,-4),
\quad \text{ and } \quad
\adeg(\alpha_2\otimes\beta_1\otimes\beta_1)=(-7,-2,-5).
\]
This implies that $c_{1i}=0$.  Similarly, we have $c_{2i}=c_{3i}=0$.
\end{proof}

In contrast to Lemma \ref{yylem6.1}, we can prove that all of the
coefficients $c_{ih}$ are zero. 
Using \SI{5}, we can find all of the $x$'s:
\[
\begin{matrix}
x_{11211} = -g_{1}^3g_{2}^3r, & 
x_{21111} = -g_{1}^4g_{2}^2, \\
x_{11112} = -g_{1}^3g_{2}^3p,& 
x_{11121} = -g_{1}^3g_{2}^3q, \\
x_{12121} = -g_{1}^2g_{2}^3r, &
x_{12211} = -g_{1}^3g_{2}^2, \\
x_{22111} = -g_{1}^3g_{2}^2p,& 
x_{12112} = -g_{1}^2g_{2}^3q, \\
x_{13121} = -g_{1}^2g_{2}^2, &
x_{13211} = -g_{1}^2g_{2}^2p,\\ 
x_{23111} = -g_{1}^2g_{2}^2q, &
x_{13112} = -g_{1}g_{2}^3r,\\
x_{14112} = -g_{1}g_{2}^2, &
x_{14211} = -g_{1}g_{2}^2q, \\
x_{14121} = -g_{1}g_{2}^2p, &
x_{24111} = -g_{1}g_{2}^2r,\\
x_{isjkh}=0& \text{for all other $(i,s,j,k,h)$.}
\end{matrix}
\]
Now we input the solutions of $\{a_{ijk},b_{stmn},c_{sh},y_{ijkh},
x_{isjkh}\}$ (which are polynomial in $g_1,g_2,p,q,r,v,w$) into \SI{6}, and 
obtain $2^7$ equations involving the variables $g_1,g_2,p,q,r,v,w$. 
This list of equations from \SI{6} is easy to generate by using 
the mathematical software Maple 8, and most of them are easy to analyze. 
We omit this list here because of its length. Instead 
we just give the non-trivial equations, which we call now the
``solutions.''
 
There are several different solutions:

\begin{solution}\label{sol1.1}
Assume $vq\neq0$. Let $f=(g_1g_2)^{-1}$. Then $v$
is a free variable and the others are functions of $v$ and $f$:
\begin{gather*}
g_1=-f^3, \quad g_2=-f^{-4}, \quad v=v, \quad w = f^2, \\
p = -f+v,\quad q = -(-f+v)f,\quad  r = -f^3.
\end{gather*}
\end{solution}

\begin{solution}\label{sol1.2}
Assume that $v=0$. There are three cases.
\begin{gather*}
g_{1} = -p^3, \quad g_{2} = -1/p^4,\quad v=0\quad 
w = -p^2, 
\tag{a}
\\
p = p, \quad q = p^2, \quad r = p^3. \\
g_{1} = i p^3, \quad
g_{2} = -1/p^4,\quad
v=0, \quad w = i p^2,
\tag{b}
\\
p = p, \quad q = p^2, \quad r = p^3,
\intertext{where $i^2=-1$.}
g_{1} = p^3, \quad 
g_{2} = -1/p^4, \quad
v = 0, \quad w = p^2, 
\tag{c}
\\
p = p, \quad
q = p^2, \quad r = p^3.
\end{gather*}
Solution~\ref{sol1.2}(c) is a special case of \ref{sol1.1} when $v=0$
and $p=-f$, so we don't need to consider this one.
\end{solution}

\begin{solution}\label{sol1.3}
Assume that $q=0$. There are two cases.
\begin{gather*}
g_1= -v^3, \quad g_2= jv^{-4},
\quad v = v,\quad w = v^2, 
\tag{a}
\\
p = 0, \quad q = 0,\quad 
r = j^{-1}v^3.
\\
\intertext{where $j^2-j+1=0$.}
g_1 = -v^3, \quad g_2 = -1/v^4,
\quad v=v, \quad w = v^2,
\tag{b}
\\
p = 0,\quad q = 0, \quad r = -v^3.
\end{gather*}
Solution~\ref{sol1.3}(b) is a special case of \ref{sol1.1} when $f=v$,
so we don't need to consider this one.
\end{solution}

We now recall the definitions of Auslander regular and Cohen-Macaulay.
For any $A$-module $M$ we define the \emph{grade} (or \emph{$j$-number}) of $M$ 
with respect to $A$ to be 
\[
j(M)=\inf\{q\;|\; \Ext^q_A(M,A)\neq 0\}.
\]
The grade of a right $A$-module is defined similarly.
An algebra $A$ is called {\it Cohen-Macaulay} if
\[
j(M)+\GKdim M=\GKdim A
\]
for all finitely generated left and right $A$-modules $M$. 
An algebra $A$ is called {\it Auslander regular} if the following
conditions hold.
\begin{enumerate}
\item
$A$ is noetherian and has finite global dimension;
\item
for every finitely generated left $A$-module $M$, every integer
$q$, and a right $A$-submodule $N\subset \Ext^q_A(M,A)$, 
one has $j(N)\geq q$;
\item
the same holds after exchanging left with right.
\end{enumerate}

Let $A$ be an algebra with automorphism $\sigma$ and 
$\sigma$-derivation $\delta$. The \emph{Ore extension}
$A[t;\sigma,\delta]$ is defined to be the algebra generated by
$A$ and $t$ subject to the following relation:
\[
at=t\sigma(a) +\delta(a)
\]
for all $a\in A$.

Now we consider Solution~\ref{sol1.1}, which includes
Solutions~\ref{sol1.2}(c) and \ref{sol1.3}(b).  Let $D(v,p)$ be the
algebra $k\langle z_1,z_2\rangle/ (r_3,r_4)$ where
\begin{gather*}
r_3=z_1z_2^2+vz_2z_1z_2+p^2 z_2^2z_1, \\
r_4=z_1^3z_2+(v+p)z_1^2z_2z_1+(p^2+pv)z_1z_2z_1^2+p^3 z_2z_1^3.
\end{gather*}
This is the algebra coming from Solution~\ref{sol1.1} by replacing the
letter $f$ in Solution~\ref{sol1.1} by $-p$ (in order to match the form in 
part (d) of Theorem \ref{yythm0.1}).

\begin{proposition}
\label{yyprop7.1} Let $D:=D(v,p)$ be defined as above. Then 
$D$ is a noetherian AS regular algebra of global dimension 4 
if and only if $p\neq 0$. Further, if $p\neq 0$, then $D$ is 
Auslander regular and Cohen-Macaulay. 
\end{proposition}

\begin{proof} 
If $p=0$, then $D$ is not a domain since $r_3$ is a product of two
lower degree polynomials.  Hence $D$ can not be AS regular of
dimension 4 by Lemma~\ref{yylem1.3}.

In the rest of the proof we assume that $p\neq 0$. Let $c$ and 
$d$ be two scalars such that $c+d=v$ 
and $cd=p^2$. We claim that $D\cong B[z_1;\sigma,\delta]$ where
\begin{gather*}
B=k\langle x,y,z_2\rangle/(yx+pxy, xz_2+d z_2x, yz_2-p^2 z_2y),
\\
\sigma: x\mapsto -d^{-1} x,\quad y\mapsto -p^{-1}y,\quad 
z_2\mapsto -c^{-1}z_2,
\\
\delta: x\mapsto d^{-1}y,\quad y\mapsto 0,\quad z_2\mapsto c^{-1}x.
\end{gather*}
Verifications of $\sigma$ being an automorphism and $\delta$ being
a $\sigma$-derivation are straightforward. Since $B$ is a skew 
polynomial ring (and an iterated Ore extension), it is AS,
Auslander regular, Cohen-Macaulay, and has global dimension 3.
Hence $R:=B[z_1,\sigma,\delta]$ is a connected graded AS regular 
algebra of global dimension 4; and  further by \cite[Theorem 4.2]{Ek} 
it is Auslander regular and Cohen-Macaulay. This also follows from the 
fact that $R$ has a filtration such that the associated graded algebra 
of $R$ is isomorphic to the skew polynomial ring 
$k_{p_{ij}}[x_1,x_2,x_3,x_4]$ for some $0\neq p_{ij}\in k$.

We may re-write the relations between $z_1$ and $x,y,z_2$ as
follows:
\[
z_1 x=(-dx)z_1+y,\quad z_1y=(-py)z_1,\quad z_1 z_2=(-cz_2)z_1+x.
\]
Then $x$ and $y$ can be generated from $z_1$ and $z_2$ as
\begin{equation}
\label{E7.1.1}
x=z_1z_2+cz_2z_1,\quad y=z_1x+dxz_1.
\end{equation}
Hence $R$ is generated by $z_1$ and $z_2$. The relation $r_3$ 
is equivalent to $xz_2+dz_2x=0$, and the relation $r_4$ is 
equivalent to $z_1y=-pyz_1$. Hence there is a surjective map 
from $D\to R$. To see that this map is an isomorphism, we define 
$x$ and $y$ in $D$ using the formula given in \eqref{E7.1.1}.
The relation $yz_2-p^2z_2y=0$ follows from $r_3$, and the 
relation $yx+pxy=0$ follows from $yz_2-p^2z_2y=0$ and $r_4$. 
Hence the defining relations of $R$ hold in $D$. 
Therefore $D$ is isomorphic to $R$.
\end{proof}

To prove that the algebras in Solution~\ref{sol1.2}(a) are AS regular, we 
need some lemmas.

\begin{lemma}
\label{yylem7.2}
Let $A$ be a connected graded algebra and $h$ a 
homogeneous normal non-zero-divisor in $A$. If $A/(h)$ has
global dimension $n$, then $A$ has global dimension $n+1$.
\end{lemma} 

\begin{proof} 
As an $A$-module, $A/(h)$ has projective dimension 1, so since $k$ has
finite projective dimension over $A/(h)$, it also has finite
projective dimension over $A$ (see \cite[8.2.2(i)]{MR}). Combined
with the fact that $A$ is connected graded, this means that
$A$ has finite global dimension.  By the Rees lemma the injective 
dimension of $A/(h)$ is one less than the injective dimension 
of $A$.  The assertion follows from the fact that the injective
dimension of $A$ is equal to the global dimension of $A$ when 
the latter is finite. 
\end{proof}

Let $O$ be the connected graded algebra $k\langle z_1,z_2\rangle/
(z_2^2z_1,z_2z_1^3, z_2z_1z_2z_1^2)$. 

\begin{lemma}
\label{yylem7.3}
The connected graded algebra $O$ has the following Hilbert series:
\[
H_O(t)=\frac{1}{(1-t)^2(1-t^2)(1-t^3)}.
\]
\end{lemma}

\begin{proof} By a result of Anick \cite[1.4]{An}, the trivial 
module $k$ has a free resolution
\[
\dotsb \to kV^{(4)}\otimes O\to kV^{(3)}\otimes O\to
kV^{(2)}\otimes O\to kV^{(1)}\otimes O\to (kz_1+kz_2)\otimes O\to k\to 0,
\]
where $V^{(n)}$ is the set of $n$-chains defined in \cite[p.\ 643]{An}. 
In this case we have  
$V^{(1)}=\{z_2^2z_1,z_2z_1^3,z_2z_1z_2z_1^2\}$,
$V^{(2)}=\{z_2^2z_1z_2z_1^2,z_2^2z_1^3,z_2z_1z_2z_1^3\}$,
$V^{(3)}=\{z_2^2z_1z_2z_1^3\}$, and $V^{(n)}=\emptyset$ for $n>3$. 
So the global dimension of $O$ is 3 and the Hilbert series of
$O$ is the inverse of the following polynomial:
\[
1-2t+(t^3+t^4+t^5)-(t^5+2t^6)+t^7=(1-t)^2(1-t^2)(1-t^3).
\] 
\end{proof}

Since $O$ is a monomial algebra, a standard $k$-linear basis of $O$ is 
the set of monomials that do not contain $z_2^2z_1$, $z_2z_1^3$, 
or $z_2z_1z_2z_1^2$ as a submonomial. Using this fact one can also 
compute the Hilbert series of $O$ by counting the monomial 
basis in each degree.

In the rest of this section we call a monomial in two
variables $z_1$ and $z_2$ {\it standard} if it 
does not contains a submonomial of the form $z_2^2z_1$, 
$z_2z_1^3$, or $z_2z_1z_2z_1^2$.  Otherwise a monomial is called 
{\it non-standard}. The set of standard monomials forms a $k$-linear
basis of $O$. 

Let $A(p)$ be the algebra $k\langle z_1,z_2\rangle/(r_3,r_4)$,
where 
\begin{gather*}
r_3=z_1z_2^2-p^2 z_2^2z_1, \\
r_4=z_1^3z_2+pz_1^2z_2z_1+p^2z_1z_2z_1^2+p^3 z_2z_1^3.
\end{gather*}
If $p=0$, then $A(p)$ is not a domain, so it is not a noetherian
AS regular algebra of dimension 4. 

\begin{lemma}
\label{yylem7.4} Let $A:=A(p)$ be defined as above. If $p\neq 0$, 
then  
\[
H_A(t)=\frac{1}{(1-t)^2(1-t^2)(1-t^3)}.
\]
\end{lemma}

\begin{proof} We follows the procedure in Bergman's diamond
lemma \cite{Be}. 

Let $z_1 < z_2$ and order all monomials using lexicographic order.
We will find a reduction system (a set of relations) such that
all ambiguities are resolvable \cite[1.2]{Be}. The relations
$r_3$ and $r_4$ can be written as
\begin{gather*}
z_2^2 z_1=p^{-2} z_1 z_2^2, \\
z_2z_1^3=-p^{-3} z_1^3z_2-p^{-2} z_1^2z_2z_1-p^{-1} z_1z_2z_1^2.
\end{gather*}
Using either $r_3$ or $r_4$, the monomial $z_2^2 z_1^3=(z_2^2z_1)z_1^2=
z_2(z_2z_1^3)$ can be reduced in two different ways. In \cite[p. 181]{Be}, 
the $5$-tuple $(r_3,r_4,z_2, z_2z_1,z_1^2)$ is called an {\it overlap 
ambiguity}. For simplicity, we call the monomial $z_2^2 z_1^3$ an 
overlap ambiguity by abuse of the notation. We now reduce 
$z_2^2z_1^3$ in two ways:
\[
z_2^2z_1^3=(z_2^2z_1)z_1^2=(p^{-2}z_1z_2^2) z_1^2=\cdots =
p^{-6}z_1^3z_2^2
\]
and 
\begin{align*}
z_2^2z_1^3 &=z_2(-p^{-3}z_1^3z_2-p^{-2}z_1^2z_2z_1-p^{-1}z_1z_2z_1^2) \\
 &=-p^{-3} z_2z_1^3z_2- p^{-2}z_2z_1^2z_2z_1-p^{-1}z_2z_1z_2z_1^2 \\
 &=-p^{-3} (-p^{-3}z_1^3z_2-p^{-2}z_1^2z_2z_1-p^{-1}z_1z_2z_1^2)z_2
- p^{-2}z_2z_1^2z_2z_1-p^{-1}z_2z_1z_2z_1^2 \\
 &=p^{-6} z_1^3z_2^2+p^{-5}z_1^2z_2z_1z_2+p^{-4}z_1z_2z_1^2z_2
- p^{-2}z_2z_1^2z_2z_1-p^{-1}z_2z_1z_2z_1^2.
\end{align*}
Using these we obtain another relation of $A$ of degree 5, called $r_5$:
\[
z_2z_1z_2z_1^2=-p^{-1} z_2z_1^2z_2z_1+p^{-3} z_1z_2z_1^2z_2+p^{-4}
z_1^2z_2z_1z_2.
\]
Starting from $r_3$, $r_4$, and $r_5$, there are three overlap ambiguities:
\[
\{z_2^2z_1^3,\quad z_2z_1z_2z_1^3, \quad z_2^2z_1z_2z_1^2\}.
\]
It is easy to verify that all ambiguities are resolvable,
which means that the two different ways of reducing the monomials
give the same result.  If $S$ is the set of non-standard
monomials, then all ambiguities of $S$ are resolvable. 
By Bergman's diamond lemma \cite[1.2]{Be}, the set of 
standard monomials forms a $k$-linear basis of $A$. Therefore the 
Hilbert series of $A$ is equal to the Hilbert series of the algebra 
$O$, which is the desired Hilbert series.
\end{proof}

\begin{proposition}
\label{yyprop7.5} Let $A:=A(p)$ be as in Lemma \ref{yylem7.4} with
$p\neq 0$. Then $A$ is a noetherian AS regular algebra of global 
dimension 4. Further, $A$ is Auslander regular and Cohen-Macaulay. 
\end{proposition}

\begin{proof} We first claim that $h=z_1^2z_2+p^2 z_2z_1^2$ is a 
normal element of $A$. The relation $r_4$ is equivalent to the
relation $z_1h+phz_1=0$. Using $r_3$, we have
\[
hz_2=z_1^2z_2^2+p^2 z_2z_1^2z_2=p^4z_2^2z_1^2+p^2z_2z_1^2z_2=
p^2 z_2h.
\]
Hence $h$ is a normal element of $A$ and 
\[
AhA=Ah=hA.
\]
Let $B=A/(h)=k\langle z_1,z_2\rangle/(r_3,r_4,h)$. Since $r_4$ is
generated by $h$ and $z_1$, $B$ is isomorphic to $B':=
k\langle z_1,z_2\rangle/(r_3,h)$. By \cite{ASc, ATV1, ATV2}, $B'$ 
is a noetherian AS regular algebra of global dimension 3 of 
type $S_2$. Thus 
\[
H_B(t)=\frac{1}{(1-t)^2(1-t^2)}=(1-t^3)H_A(t).
\]
Since $B=A/hA=A/Ah$, we have $H_{hA}(t)=H_{Ah}(t)=t^3 H_A(t)$.
This implies that $h$ is a non-zero-divisor. By Lemma \ref{yylem7.2},
$A$ has global dimension 4. By \cite[8.2]{ATV1}, $A$ is 
noetherian. Since $B=A/(h)$ is Auslander regular
and Cohen-Macaulay, by \cite[5.10]{Lev}, so is $A$.
\end{proof}

We now consider the algebras in Solution~\ref{sol1.2}(b).  Let $B(p)$
be the algebra $k\langle z_1,z_2\rangle/(r_3,r_4)$, where
\begin{gather*}
r_3=z_1z_2^2+i p^2 z_2^2z_1, \\
r_4=z_1^3z_2+pz_1^2z_2z_1
+p^2z_1z_2z_1^2+p^3 z_2z_1^3,
\end{gather*}
and where $i^2=-1$.
The difference between the algebra $A(p)$ in Lemma \ref{yylem7.4} 
and $B(p)$ is the relation $r_3$. Suppose $p\neq 0$. Let $S$ be 
the set of standard monomials.
We list some basic properties of $B(p)$ below.

\begin{lemma}
\label{yylem7.6} Let $B:=B(p)$ be defined as above and let $p\neq 0$. 
\begin{enumerate}
\item 
$z_2^2$ is a normal element of $B$. $z_1^4$ is a normal element in $B$ and
is central if $p^4=1$.
\item
$B/(z_2^2)$ is a factor ring of the AS 
regular algebra in Lemma \ref{yylem7.4}. 
$B$ is noetherian.
\item 
The automorphism group of $B$ is $(k-\{0\})^{\times 2}$, acting separately
on $z_1$ and $z_2$.
\item
$B(p)$ is a twist of $B(1)$.
\item
The Hilbert series of $B$ is $((1-t)^2(1-t^2)(1-t^3))^{-1}$.
The set $S$ is a $k$-linear basis of $B$.
\item
$z_1$ and $z_2$ are non-zero-divisors of $B$. 
\item
$z_1$ is a non-zero-divisor in the factor ring $B/(z_2^2)$.
At the level of $B$, if $cz_1\in Bz_2^2$, then $c\in Bz_2^2$.
\item
$B$ has global dimension 4.
\end{enumerate}
\end{lemma}

\begin{proof} (a) It follows from the relations $r_3$ and $r_4$
respectively that $z_1z_2^2=-i p^2 z_2^2z_1$ and 
$z_1^4 z_2=p^4 z_2 z_1^4$.

(b) Let $A$ be the noetherian algebra in Lemma \ref{yylem7.4}.  Then 
$B/(z_2^2)\cong A/(z_2^2)$, and this quotient is noetherian.  By
\cite[8.2]{ATV1} $B$ is noetherian.

(c) Since the relations $r_3,r_4$ are ${\mathbb Z}^2$-homogeneous,
the ``diagonal'' map 
\[
z_1\mapsto a_{11}z_1, \quad z_2\mapsto a_{22}z_2
\] 
is an automorphism of $B$. It remains to show every automorphism 
is of this form. For any automorphism $\sigma$, $\sigma(r_3)$ must 
be a scalar multiple of $r_3$ and $\sigma(r_4)$ must be
generated by $r_3$ and $r_4$. 

Suppose $\sigma: z_1\mapsto a_{11}z_1+a_{12}z_2,\quad
z_2\mapsto a_{21}z_1+a_{22}z_2$ is an automorphism of $A$. We write
$\sigma(r_3)$ as a linear combination of monomials:
\[
\sigma(r_3)=a_{11}a_{21}^2(1+ip^2)z_{1}^3+\cdots + 
a_{12}a_{22}^2(1+ip^2)z_2^3.
\]
Since $\sigma(r_3)$ is a scalar multiple of $r_3$, the coefficients
of $z_1^3$ and $z_2^3$ are zero. 

If $1+ip^2=0$, then $p^4=-1$. This is one case we need to
consider later.

If $1+ip^2\neq 0$ or $p^2\neq i$, then $a_{11}a_{21}=0=a_{12}a_{22}$. 
Since $\sigma$ is an automorphism, there are two cases to consider:
either $a_{21}=a_{12}=0$ or $a_{11}=a_{22}=0$. 
The first case gives the diagonal automorphism. In the second case 
$\sigma(r_3)\neq cr_3$ for any $c\in k-\{0\}$. So $\sigma$ is not
an automorphism. 

It remains to consider the case of $p^4=-1$. Let $\Omega=
1+p+p^2+p^3$, which is nonzero when $p^4=-1$. Express $\sigma(r_4)$
as a linear combination of monomials:
\[
\sigma(r_4)=\Omega a_{11}^3a_{21}z_1^4+\cdots+ \Omega a_{12}^3a_{22} 
z_2^4=:r'.
\]
To make $r'=0$ in $A$, we need both $a_{11}^3a_{21}$ and $a_{12}^3a_{22}$
to be zero. Thus we have either $a_{21}=a_{12}=0$ or $a_{11}=a_{22}=0$.
Both were discussed in the previous paragraph.

(d) Via the automorphism $\tau: z_1\mapsto z_1, \ z_2\mapsto p z_2$, 
the twisted algebra $B(p)^\tau$ is isomorphic to $B(1)$. 
(For definitions and details see \cite{Zh1}.)

(e) These can be proved by imitating the proof of 
Lemma \ref{yylem7.4}. The computation is slightly more 
complicated since there is an extra term in $r_5$.
The proof is omitted. 

(f) If $f\in S$, then $f$ does not contain a submonomial 
of the forms $z_2^2z_1$, $z_2z_1^3$, or $z_2z_1z_2z_1^2$.
This implies that $z_1f$ and $fz_2$ do not contain a submonomial 
of the form $z_2^2z_1$, $z_2z_1^3$, or $z_2z_1z_2z_1^2$.
Hence $z_1f, fz_2\in S$. 

Let $a$ be a non-zero element in $B$. Write $a$ as a linear 
combinations of monomials in $S$, say $a=\sum c_i f_i$ 
where $c_i\in k-\{0\}$ and the $f_i$ are distinct elements in 
$S$. Then $z_1 a= \sum_i c_i z_1 f_i$ where $\{z_1 f_i\}$ 
is a set of different elements in $S$. Hence $z_1a\neq 0$. 
This says that $z_{1}$ is a left non-zero-divisor. Since $z_1^4$
is a normal by part (a), $z_1$ is also a right non-zero-divisor. 

A similar argument works for $z_2$.

(g) We use Bergman's diamond lemma \cite{Be}. Write
the two relations of $B/(z_2^2)$ as
\begin{gather*}
z_2^2=0, \\
z_2z_1^3=-p^{-3}z_1^3z_2-p^{-2}z_1^2z_2z_1-p^{-1}z_1z_2z_1^2.
\end{gather*}
As in the proof of Lemma \ref{yylem7.4}, we can reduce the 
overlap ambiguity $z_2^2z_1^3$ in two different ways and obtain 
a new relation, 
\[
z_2z_1z_2z_1^2=-p^{-1}z_2z_1^2z_2z_1+p^{-3}z_1z_2z_1^2z_2
+p^{-4} z_1^2z_2z_1z_2.
\]
Using these three relations (two original and the new one displayed
above) one can verify that all three 
overlap ambiguities $\{z_2^2z_1^3, z_2z_1z_2z_1^3,
z_2^2z_1z_2z_1^2\}$ are resolvable and therefore they will 
not create new relations. By Bergman's diamond lemma
\cite[1.2]{Be}, the set $S/(z_2^2)$ of all monomials 
that do not contain a submonomial of the form $z_2^2, z_2z_1^3, 
z_2z_1z_2z_1^2$ is a $k$-linear basis of $B/(z_2^2)$.
If $f=z_1^{i_1}z_2\cdots z_2z_1^{i_n}$ is an
element in $S/(z_2^2)$, then so is $z_1f$. So $z_1$ is a left 
non-zero-divisor. Since $z_1^4$ is normal, $z_1$ is also a 
right non-zero-divisor.

(h) Since $B(p)$ is a twisted algebra of $B(1)$, $B(p)$ and $B(1)$ have
the same global dimension \cite[5.7 and 5.11]{Zh1}. It suffices to show that 
$B(1)$ has global dimension 4. So we assume that $p=1$.
(When $p\neq 1$ one can also give a similar proof, but $p$ and the powers of $p$
appear in various places.) 

Let $R=B(1)$ and consider the following complex of left $R$-modules
\begin{equation}
\label{E7.6.1}
0\to R\delta^* \to R\gamma_1^*\oplus R\gamma_2^* \to 
R\beta_1^* \oplus R\beta_2^* \to R\alpha_1^*\oplus R\alpha_2^*
\to R\to  k\to 0
\end{equation}
where 

the map $R\to  k$ sends 
\[
z_1\mapsto 0, \quad z_2 \mapsto 0;
\]

the map $R\alpha_1^*\oplus R\alpha_2^*\to R$ sends 
\[
\alpha_1^*\mapsto z_1,\quad  \alpha_2^*\mapsto z_2;
\]

the map $R\beta_1^* \oplus R\beta_2^* \to R\alpha_1^*\oplus R\alpha_2^*$
sends 
\[
\beta_1^*\mapsto i z_2^2 \alpha_1^*+z_1z_2 \alpha_2^*,\quad 
\beta_2^*\mapsto (z_1^2z_2+z_1z_2z_1+z_2z_1^2)\alpha_1^*+z_1^3\alpha_2^*;
\]

the map $R\gamma_1^*\oplus R\gamma_2^* \to R\beta_1^* \oplus R\beta_2^* $
sends 
\[
\gamma_1^*\mapsto (z_1^2z_2+z_1z_2z_1+z_2z_1^2)\beta_1^*+
i z_2^2 \beta_2^*,\quad 
\gamma_2^*\mapsto z_1^3 \beta_1^*+(-z_2z_1)\beta_2^*;
\]

and the map $R\delta^* \to R\gamma_1^*\oplus R\gamma_2^* $ sends 
\[
\delta^*\mapsto z_1\gamma_1^*+z_2\gamma_2^*.
\] 

Using the relations $r_3$ and $r_4$, it is easy to verify that
\eqref{E7.6.1} is a complex, namely, the composition of any two
consecutive maps is zero. Since $z_1$ and $z_2$ are the generators and $r_3$
and $r_4$ are the relations, the complex \eqref{E7.6.1} is exact at the
positions $k$, $R$ and $R\alpha_1^*\oplus R\alpha_2^*$. Since $z_1$
and $z_2$ are left non-zero-divisors, the complex is exact at
$R\delta^*$.

Next we want to show that \eqref{E7.6.1} is exact at the position
$R\gamma_1^*\oplus R\gamma_2^* $. Let $F$ be the map $R\gamma_1^*\oplus 
R\gamma_2^* \to  R\beta_1^* \oplus R\beta_2^*$.  Exactness 
is equivalent to 
\[
\ker F=\{az_1\gamma_1^*+az_2\gamma_2^* \;| \; a\in R\}=:I.
\]
It is clear that $\ker F\supset I$. It remains to show that 
$\ker F\subset I$. If $f=f_1 \gamma_1^* +f_2 \gamma_2^* \in 
\ker F$, then the $\beta_2^*$ component of the equation 
$F(f)=0$ is
\[
f_1 i z_2^2 +f_2 (-z_2z_1)=0.
\]
So $(f_2z_2) z_1=if_1z_2^2\in Rz_2^2$. By (g) $f_2z_2\in Rz_2^2$.
Since $z_2$ is a non-zero-divisor, $f_2\in Rz_2$. Write $f_2=az_2$.
Using $r_3$ and the equation above, $f_1=az_1$ and $f\in I$.

Applying the additive function $H_M(t)$ to the complex 
\eqref{E7.6.1} gives
\[
H_R(t)p(t)=1-H_{H^{-2}}(t),
\]
where $p(t)=(1-t)^2(1-t^2)(1-t^3)$ and where $H^{-2}$ is the 
cohomology of the complex at the position $R\beta_1^*\oplus R\beta_2^*$.
By (e), $H_R(t)p(t)=1$. Hence $H_{H^{-2}}(t)=0$ and therefore 
$H^{-2}=0$. This implies that \eqref{E7.6.1} is exact. So $R=B(1)$ 
has global dimension 4.
\end{proof}

\begin{proposition}
\label{yyprop7.7} Let $B=B(p)$ be the algebra in Lemma \ref{yylem7.6}.
Then $B$ is Auslander regular and Cohen-Macaulay of global
dimension 4.
\end{proposition}

\begin{proof} By Lemma \ref{yylem7.6}(b), $B/(z_2^2)$ is isomorphic
to $A/(z_2^2)$, where $A=A(p)$ is the algebra in Lemma~\ref{yylem7.4}. 
Since $z_2^2$ is a normal non-zero-divisor of $B$ and of $A$, and since
$A$ is an Auslander (and AS) Gorenstein and Cohen-Macaulay, by using 
\cite[5.10]{Lev} twice, $B$ is Auslander (and AS) Gorenstein and 
Cohen-Macaulay.  The regularity follows from Lemma \ref{yylem7.6}(h).
\end{proof}

By a direct computation, the algebra $A(p)$ is not isomorphic to 
$B(p')$ for all $p,p'$. Another way of showing this is to use the 
following lemma. Recall that there is a surjective map from $A(p)$ 
to an AS regular algebra of global dimension 3 (see the proof of 
Proposition \ref{yyprop7.5}).

\begin{lemma}
\label{yylem7.8} 
Let $B=B(p)$ be the AS regular algebra in Lemma \ref{yylem7.6}. 
\begin{enumerate}
\item There are no nonzero normal elements of degree 3 in $B$.
\item There is no surjective map from $B$ to any AS regular algebra of global
dimension 3.
\end{enumerate}
\end{lemma} 

\begin{proof}
(a) Suppose that $a$ is a nonzero normal element of degree 3. Then
there is an automorphism $\sigma$ such that $ba=ab^\sigma$ for all
$b\in B$. By Lemma \ref{yylem7.6}(c) there are nonzero 
scalars $c_1$ and $c_2$ such that
\[
z_1 a=c_1 az_1 \qquad\text{and}\qquad z_2 a=c_2 az_2.
\]
Since $B$ is ${\mathbb Z}^2$-graded, we may assume $a$ is 
${\mathbb Z}^2$-homogeneous. Then up to a scalar, $a$ is one of
the following:
\[
z_1^3, \quad b_1 z_1^2z_2+b_2 z_1z_2z_1+b_3 z_2z_1^2,
\quad d_1 z_1z_2^2+d_2 z_2z_1z_2, \quad z_2^3.
\]

If $a=z_1^3$, then $z_2 a\neq c_2 az_2$ for any $c_2\in k$. 

If $a= b_1 z_1^2z_2+b_2 z_1z_2z_1+b_3 z_2z_1^2$, then 
$z_2 a=c_2 az_2$ implies that $b_2=0$, $c_2=ip^{-2}$ and 
$b_3= -i p^{-2}b_1$.
In this case $z_1 a\neq c_1 az_1$ for any $c_1\in k$.

If $a=d_1 z_1z_2^2+d_2 z_2z_1z_2 $, then $az_1=c_1 z_1a$ implies that
$d_2=0$, and $az_2 =c_2 z_2 a$ implies that $d_1=0$.

If $a= z_2^3$, then $z_1a \neq c_1 az_1$ for any $c_1\in k$.

(b) If there is a surjective map from $B$ to a 3-dimensional
AS regular algebra $D$, then $D$ is generated by two degree 1 elements.
By Artin and Schelter's classification \cite{ASc}, $H_D(t)=((1-t)^2(1-t^2))^{-1}$.
Comparing the Hilbert series of $B$ and $D$ we see that 
$D\cong B/(a)$ where $a$ is a degree 3 normal non-zero-divisor.
This contradicts (a). 
\end{proof}

The next case is Solution~\ref{sol1.3}(a).  Let $C(p)$
be the algebra $k\langle z_1,z_2\rangle/(r_3,r_4)$ with
\begin{gather*}
r_3= z_1z_2^2 +p z_2 z_1 z_2 + p^2 z_2^2 z_1, \\
r_4= z_1^3 z_2 + j p^3 z_2 z_1^3,
\end{gather*}
where $p\neq 0$ and $j^2-j+1=0$.

Similar to Lemma \ref{yylem7.6}, we can show the following.
We will omit most of details in the proof since the ideas are 
the same as in the proof of Lemma \ref{yylem7.6}. 

\begin{lemma}
\label{yylem7.9} Let $C=C(p)$ be defined as above. 
\begin{enumerate}
\item 
$z_1^3$ and $z_2^3$ are normal elements of $C$.
\item
$C/(z_1^3)$ is a factor ring of the AS regular algebra $D(-p,p)$ in 
Proposition \ref{yyprop7.1} 
(with $v=-p$). As a consequence, $C$ is noetherian.
\item 
The automorphism group of $C$ is $(k-\{0\})^{\times 2}$,
the factors acting separately on $z_1$ and $z_2$. 
\item
$C(p)$ is a twist of $C(1)$.
\item
The Hilbert series of $C$ is $((1-t)^2(1-t^2)(1-t^3))^{-1}$.
The set $S$ of standard monomials is a $k$-linear basis of $C$.
\item
$z_1$ and $z_2$ are non-zero-divisors of $C$.
\item
$z_2$ is a non-zero-divisor of $C/(z_1^3)$. At the level 
of $C$, if $cz_2\in Cz_1^3$, then $c\in Cz_1^3$. 
\item
$C$ has global dimension 4.
\end{enumerate}
\end{lemma}

\begin{proof} 
(a) It follows from $r_4$ that $z_1^3$ is normal. To see 
that $z_2^3$ is normal, we use the relation $r_3$:
\[
z_1z_2^3=(-pz_2z_1z_2-p^2 z_2^2z_1)z_2=
-pz_2[-pz_2z_1z_2-p^2 z_2^2z_1] -p^2z_2^2z_1z_2=p^3 z_2^3z_1.
\]

(b,c,d) Similar to Lemma \ref{yylem7.6}(b,c,d).

(e) We give a sketch of the proof when $p=1$. The relations
$r_3,r_4$ can be written as 
\begin{gather*}
z_2^2 z_1 =-z_1 z_2^2- z_2 z_1 z_2, \\
z_2 z_1^3=-j^{-1} z_1^3 z_2.
\end{gather*}
The first overlap ambiguity is $z_2^2 z_1^3$. The extra condition 
$r_5$ derived from this overlap ambiguity is
\[
z_2 z_1 z_2 z_1^2= j^2 z_1^3 z_2^2 -z_1^2 z_2 z_1 z_2 + z_1 z_2  
z_1 z_2 z_1.
\]
Using these three relations, all ambiguities are resolvable. Hence 
$S$ is a $k$-linear basis of $C$. The assertion follows. 

(f,g) Similar to Lemma \ref{yylem7.6}(f,g).

(h) We use the same method as in the proof of Lemma 
\ref{yylem7.6}(h). Again we assume that $p=1$. Let $R=C(1)$ and consider
the following complex. 

\begin{equation}
\label{E7.9.1}
0\to R\delta^* \to R\gamma_1^*\oplus R\gamma_2^* \to 
R\beta_1^* \oplus R\beta_2^* \to R\alpha_1^*\oplus R\alpha_2^*
\to R\to  k\to 0
\end{equation}
where 

the map $R\to  k$ sends 
\[
z_1\mapsto 0, \quad z_2\mapsto 0;
\] 

the map $R\alpha_1^*\oplus R\alpha_2^*\to R$ sends 
\[
\alpha_1^*\mapsto z_1, \quad  \alpha_2^*\mapsto z_2;
\]

the map $R\beta_1^* \oplus R\beta_2^* \to R\alpha_1^*\oplus R\alpha_2^*$
sends 
\[
\beta_1^*\mapsto z_2^2 \alpha_1^*+(z_1z_2 + z_2z_1) \alpha_2^*,
\quad \beta_2^*\mapsto j z_2 z_1^2 \alpha_1^*+z_1^3\alpha_2^*;
\]

the map $R\gamma_1^*\oplus R\gamma_2^* \to R\beta_1^* \oplus R\beta_2^* $
sends 
\[
\gamma_1^*\mapsto (-j z_1^2z_2)\beta_1^*+ z_2^2 \beta_2^*,
\quad \gamma_2^*\mapsto (-j^2 z_1^3) \beta_1^*+(z_1z_2+z_2z_1)\beta_2^*;
\]

and the map $R\delta^* \to R\gamma_1^*\oplus R\gamma_2^* $ sends 
\[
\delta^*\mapsto z_1\gamma_1^*+z_2\gamma_2^*.
\] 

It is straightforward to check that this is a complex and 
that it is exact at positions $k$, $R$, $R\alpha_1^*\oplus R\alpha_2^*$
and $R\delta^*$. Next we check that it is exact at position 
$R\gamma_1^*\oplus R\gamma_2^*$. This is equivalent to
checking that $\ker F= \{a z_1\gamma_1^*+a z_2\gamma_2^*| a\in C\}$,
where $F$ is the map from $R\gamma_1^*\oplus R\gamma_2^* 
\to R\beta_1^* \oplus R\beta_2^* $. 
Let $f=f_1 \gamma_1^* +f_2 \gamma_2^* \in \ker F$. Then 
\[
0=F(f)=(-j)[f_1 z_1^2 z_2 + f_2 j z_1^3]\beta_1^*+
[f_1 z_2^2 +f_2 (z_1z_2+z_2z_1)]\beta_2^*.
\]
Hence $f_1 z_1^2 z_2 + f_2 j z_1^3=0$ and $(f_1z_1^2) z_2\in
Rz_1^3$. By (g), $(f_1z_1^2)\in Rz_1^3$ and this implies that
$f_1=az_1$ for some $a\in R$. Using $f_1 z_1^2 z_2 + f_2 j z_1^3=0$
and $r_4$, we see that $f_2=az_2$ as desired.

The rest of the proof is the same as the proof of 
Lemma \ref{yylem7.6}(h).
\end{proof}

Similar to Proposition \ref{yyprop7.7}, we have the following.

\begin{proposition}
\label{yyprop7.10} Let $C=C(p)$ be the algebra in Lemma \ref{yylem7.9}.
Then $C$ is Auslander regular and Cohen-Macaulay of global
dimension 4.
\end{proposition}

\section{Properties of generic algebras}
\label{yysect7}

In this section, we prove our main theorems.


For each AS regular algebra $A$ discussed in the previous section --
$A(p)$, $B(p)$, $C(p)$, and $D(v,p)$ -- one can study the geometric
properties and invariants of its projective spectrum $\Proj A$. For
example, one can try to compute the point-scheme and the
line-scheme. To ensure that the point-scheme and the line-scheme are
in fact projective schemes it is sufficient to prove that the algebra
$A$ is {\it strongly noetherian}, namely, for every commutative
noetherian $k$-algebra $R$, $A\otimes R$ is noetherian
\cite[p.580]{ASZ}.  This property can be proved.

\begin{proposition}
\label{yyprop7.11}
Algebras $A(p)$, $B(p)$, $C(p)$ and  $D(v,p)$ are strongly noetherian.
\end{proposition}

\begin{proof} Since $k$ is trivially strongly noetherian and
$D(v,p)$ is an iterated Ore extension starting from $k$, 
by \cite[4.1(1)]{ASZ}, $D(v,p)$ is strongly noetherian.

By the proof of Proposition \ref{yyprop7.5}, $A':=A(p)/(h)$ is an AS regular 
algebra of dimension 3. Using the relations of $A'$ one sees
that $z_1^2$ and $z_2^2$ are normal elements of $A'$, and
that $A'/(z_1^2,z_2^2)=k\langle z_1,z_2\rangle/(z_1^2,z_2^2)$
is a noetherian affine PI ring. It follows from \cite[4.9(1,5)]{ASZ}
that $A(p)$ is strongly noetherian. Similar ideas apply
to $B(p)$ and $C(p)$ by using Lemmas \ref{yylem7.6}(b) and 
\ref{yylem7.9}(b).
\end{proof}

\begin{proof}[Proof of Theorems \ref{yythm0.1} and \ref{yythm0.3}]
By Propositions \ref{yyprop7.5}, \ref{yyprop7.7}, 
\ref{yyprop7.10} and \ref{yyprop7.1}, the algebras $A(p)$,
$B(p)$, $C(p)$ and $D(v,p)$ are AS regular of global dimension 
four and Auslander regular and Cohen-Macaulay.  By Proposition
\ref{yyprop7.11} these algebras are strongly noetherian.

Now let $k$ be an algebraically closed field.  By
Proposition~\ref{prop5.new} and the discussion in
Section~\ref{yysect-reg}, these are the only possible
$(m_2,m_3)$-generic, AS regular algebras of type \type{12221}.
\end{proof}

\begin{proposition}
\label{yyprop7.12} Two algebras in Theorem \ref{yythm0.1} are
isomorphic if and only if their relations are the same.
\end{proposition}

\begin{proof} Let $A=k\langle z_1,z_2\rangle/(r_3,r_4)$ 
be an algebra of the form given in Theorem \ref{yythm0.1}.
The relations $r_3$ and $r_4$ are the two relations of degree 3 and
degree 4 given in Theorem \ref{yythm0.1}. Let 
$A'=k\langle z_1,z_2\rangle/(r'_3,r'_4)$ be another such 
algebra and let $\sigma: A'\to A$ be an isomorphism of graded
algebras. By an argument similar to the one given in the proof 
of Lemmas \ref{yylem7.6}(c) and \ref{yylem7.9}(c), $\sigma(z_1)
=a_{11} z_1$ and $\sigma(z_2)=a_{22}z_2$. Using the special 
form of the relations $r_3$ and $r_4$, it easily follows that 
$r_3=r'_3$ and $r_4=r'_4$. 
\end{proof}

\begin{proof}[Proof of Theorem \ref{yythm0.2}] 
Let $A$ be a ${\mathbb Z}^2$-graded AS regular algebra of 
type \type{12221} generated by elements $z_{1}$ and $z_{2}$ in degrees
$(1,0)$ and $(0,1)$, respectively, and let $E$ be its $\Ext$-algebra.
Then we may choose a basis $\{\alpha_1,\alpha_2\}$ of $E^1$ so that 
$\adeg \alpha_1=(-1,0)$ and $\adeg \alpha_2=(0,-1)$.

By Proposition \ref{yyprop1.4}(a), $A$ has two relations 
$r_3$ and $r_4$ of Adams degree 3 and 4 respectively.
Since $A$ is a domain, $r_3$ is either $z_1z_2^2
+v z_2z_1z_2+w z_2^2 z_1$ (with $w\neq 0$) or $z_2z_1^2
+v z_1z_2z_1+w z_1^2 z_2$. By symmetry we may assume that
\[
r_3=z_1z_2^2+v z_2z_1z_2+w z_2^2 z_1.
\]
To determine $r_4$ we need to use information from the Ext-algebra.
Since $A$ is AS regular, $E$ is Frobenius (Theorem~\ref{yythm1.8}).
Let $\delta$ be a nonzero element in $E^{4}$, and let $\{\gamma_{1},
\gamma_{2} \}$ be a homogeneous basis for $E^{3}_{-6}$.  Using the
${\mathbb Z}^2$-grading, we may choose these so that $\alpha_i
\gamma_j=\gamma_j\alpha_i=0$ when $i\neq j$.  Choosing $\gamma_i$
properly, we may assume that $\alpha_i\gamma_i =\delta$ and
$\gamma_i\alpha_i=g_i \delta$.  This means that the $\Lambda$-matrix
of $E$ is always diagonal.

Since $r_4$ is homogeneous with respect to the ${\mathbb Z}^2$-graded,
$\adeg r_4$ is equal to one of the following:

\begin{center}
Case 1: $(4,0)$, $\ $
Case 2: $(0,4)$, $\ $
Case 3: $(1,3)$, $\ $
Case 4: $(2,2)$, $\ $
Case 5: $(3,1)$.
\end{center}

\noindent
These correspond exactly to the cases listed in
Proposition~\ref{prop5.new}.  Cases 1 and 2
contradict \eqref{E4.R4}. Case 3 contradicts \eqref{E4.R4} 
by our choice of \eqref{E5.2.1}. 
The discussion of Cases 4 and 5 gives the list of possible $E$s. 
Only those in Case 5 are possible $E$s for AS regular algebras of 
type \type{12221}. 
\end{proof}

Finally we mention the following proposition without proof.
The proof we have is computational. 

\begin{proposition}
\label{yyprop7.13} 
The $A_\infty$-Ext-algebras of the AS regular
algebras $D(v,p)$, $A(p)$, $B(p)$, $C(p)$ are the ones 
described in Solutions~\ref{sol1.1}, \ref{sol1.2}(a),
\ref{sol1.2}(b), and \ref{sol1.3}(a), respectively. 
\end{proposition}

\appendix
\section{$A_\infty$-structure on Ext-algebras}
\label{app}

The main goal of this appendix is to prove Theorem \ref{yythm2.2} and
Proposition~\ref{yyprop2.3}.  We also review some material which may
be useful for people working in graded ring theory.

\subsection{Kadeishvili's theorem and Merkulov's construction}
\label{sectapp1}

Let $A$ and $B$ be two $A_\infty$-algebras. A {\it morphism} of 
$A_\infty$-algebras $f: A\to B$ is a family of $k$-linear graded maps
\[
f_n: A^{\otimes n}\to B
\]
of degree $1-n$ satisfying the following {\it Stasheff morphism identities}:
\[
\sum (-1)^{r+st} f_u(id^{\otimes r}\otimes m_s\otimes 
id^{\otimes t})=\sum (-1)^{w} m_q(f_{i_1}\otimes f_{i_2}
\otimes \cdots \otimes f_{i_q})
\tag*{\MI{n}}
\]
for all $n\geq 1$, 
where the first sum runs over all decompositions $n=r+s+t$ with 
$s\geq 1$ and $r,t\geq 0$, where $u=r+1+t$, and the second sum runs over 
all $1\leq q\leq n$ and all decompositions $n=i_1+\cdots + i_{q}$
with all $i_s\geq 1$.  The sign on the right-hand side is given by
\[
w=(q-1)(i_1 -1) + (q-2) (i_2-1)+ \cdots + 2(i_{q-2}-1)+ (i_{q-1}-1).
\]

When $A$ and $B$ have a strict unit (as we always assume), 
an $A_\infty$-morphism is also required to satisfy the following extra 
{\it unital morphism conditions}:
\[
f_1(1_A)=1_B
\]
where $1_A$ and $1_B$ are strict units of $A$ and $B$ respectively, and
\[
f_n(a_1 \otimes \cdots \otimes a_n)= 0 
\]
if $n\geq 2$ and $a_i=1_A$ for some $i$.

If $A$ and $B$ have Adams gradings indexed by the same group, then 
the maps $f_i$ are required to preserve the Adams degree.

A morphism $f$ is called a {\it quasi-isomorphism} if $f_1$ is a 
quasi-isomorphism.  A morphism is {\it strict} if $f_i=0$ for all $i\neq 1$. 
The \emph{identity morphism} is the strict morphism $f$ such that $f_1$ is
the identity of $A$. When $f$ is a strict morphism from $A$ to $B$,
then the identity \MI{n} becomes
\[
f_1 m_n=m_n(f_1\otimes \cdots \otimes f_1).
\]
A morphism $f=(f_i)$ is called a {\it strict isomorphism} if it is
strict with $f_1$ a vector space isomorphism.

Let $A$ be an $A_\infty$-algebra.  Its cohomology ring is defined to be
\[
HA:=\ker m_1/\im m_1.
\] 
The following result, due to Kadeishvili \cite{Ka}, is a basic and
important property of $A_\infty$-algebras.

\begin{theorem}
\label{yythmA.1} \cite{Ka}
Let $A$ be an $A_\infty$-algebra and let $HA$ be the cohomology ring of 
$A$. There is an $A_\infty$-algebra structure on $HA$ with $m_1=0$,
constructed from the $A_\infty$-structure of $A$, such that there is a 
quasi-isomorphism of $A_\infty$-algebras $HA\to A$ lifting the identity
of $HA$.  This $A_\infty$-algebra structure on $HA$ is unique up to
quasi-isomorphism.
\end{theorem}




Kadeishvili's construction is very general. We would like to describe
some specific $A_\infty$-structures that we can work with.  Merkulov
constructed a special class of higher multiplications for $HA$ in
\cite{Me1}, in which the higher multiplications can be defined
inductively; this way, the $A_\infty$-structure can be described
more explicitly, and hence used more effectively.  For our purposes we
will describe a special case of Merkulov's construction, assuming that
$A$ is a DGA.

Let $A$ be a DGA with differential $\partial$ and multiplication $\cdot$. 
Denote by $B^n$ and $Z^n$
the coboundaries and cocycles of $A^n$, respectively. Then there
are subspaces $H^n$ and $L^n$ such that 
\[
Z^n=B^n\oplus H^n
\]
and 
\begin{equation}
\label{A.1.1}A^n=Z^n\oplus L^n=B^n\oplus H^n\oplus L^n.
\end{equation}
We will identify $HA$ with $\boplus_{n} H^n$, or embed $HA$ into $A$ by 
cocycle-sections $H^n\subset A^n$. {\it There are many different 
choices of $H^n$ and $L^n$.} 

Note that if $A$ has an Adams grading, then the decompositions above
will be chosen to respect the Adams grading, and all
maps constructed below will preserve the Adams grading.

Let $p=Pr_H: A\to A$ be a projection to $H:=\boplus_n H^n$, and 
let $Q: A\to A$ be a homotopy from $id_A$ to $p$. Hence 
we have $id_A-p=\partial Q+Q\partial$.  The map $Q$ is not unique. From now on 
we choose $Q$ as follows: 
for every $n$, 
$Q^n: A^n\to A^{n-1}$ is defined by
\begin{itemize}
\item 
$Q^n=0$ when restricted to $L^n$ and $H^n$, and 
\item
$Q^n=(\partial^{n-1}|_{L^{n-1}})^{-1}$ when restricted to $B^n$.
\end{itemize}
So the image of $Q^n$ is $L^{n-1}$. 
It follows that $Q^{n+1}\partial^n=Pr_{L^n}$ and 
$\partial^{n-1}Q^n=Pr_{B^n}$.

Define a sequence of linear maps $\lambda_n: A^{\otimes n}\to A$ of
degree $2-n$ as follows.  There is no map $\lambda_{1}$, but we
formally define the ``composite'' $Q \lambda_{1}$ by $Q \lambda_{1} =
-id_{A}$.  $\lambda_2$ is the multiplication of $A$, namely, 
$\lambda_2(a_1\otimes a_2)=a_1\cdot a_2$.
For $n\geq 3$,
$\lambda_n$ is defined by the recursive formula
\begin{equation}
\label{A.1.2}
\lambda_n=\sum_{\substack{s+t=n, \\ s,t\geq 1}} (-1)^{s+1} \lambda_2[Q\lambda_s\otimes
Q\lambda_t].
\end{equation}

We abuse notation slightly, and use $p$ to denote both the map $A
\rightarrow A$ and also (since the image of $p$ is $HA$) the map
$A \rightarrow HA$; we also use $\lambda_{i}$ both for the map
$A^{\otimes i} \rightarrow A$ and for its restriction $(HA)^{\otimes
i} \rightarrow A$.

Merkulov reproved Kadeishvili's result in \cite{Me1}.

\begin{theorem}
\label{yythmA.2} \cite{Me1}
Let $m_i=p\lambda_i$. Then $(HA, m_2,m_3, \cdots)$ is an $A_\infty$-algebra.
\end{theorem}

We can also display the quasi-isomorphism between $HA$ and $A$ directly.

\begin{proposition}
\label{yypropA.3} 
Let $\{\lambda_n\}$ be defined as above. Let $f_i=-Q 
\lambda_i: (HA)^{\otimes i}\to A$ and let $m_i=p \lambda_i:
(HA)^{\otimes i}\to HA$. Then $(HA, m_2, m_3, \cdots)$ is an 
$A_\infty$-algebra
and $f:=\{f_i\}$ is a quasi-isomorphism of $A_\infty$-algebras.
\end{proposition}

\begin{proof} 
This construction of $\{m_i\}$ and $\{f_i\}$ is a special case of
Kadeishvili's construction.
\end{proof}

Any $A_\infty$-algebra constructed as in Theorem \ref{yythmA.2} and 
Proposition \ref{yypropA.3} 
is called a {\it Merkulov model} of $A$, denoted by $H_{Mer}A$. 
The particular model depends on the decomposition \eqref{A.1.1}, but 
all Merkulov models of $A$ are quasi-isomorphic to each other. 
If $A$ has an Adams grading, then by construction all maps $m_i$ and $f_i$ 
preserve the Adams degree. 

Next we consider the unital condition. 

\begin{lemma}
\label{yylemA.4} Suppose $H^0$ is chosen to contain the unit element of $A$. 
Then $H_{Mer}A$ satisfies the strictly unital condition, and the morphism 
$f=\{f_i\}$ satisfies the unital morphism conditions.
\end{lemma}

\begin{proof} First of all, $1\in H^0$ is a unit with respect 
to $m_2$.  We use induction on $n$ to show the following, for $n \geq 3$:

$(a)_n$: $\ $ $f_{n-1}(a_1 \otimes \dots \otimes a_{n-1})=0$ if $a_i=1$ 
for some $i$.

$(b)_n$: $\ $ $\lambda_n(a_1 \otimes \dots \otimes a_n)\in L:=\boplus_n L^n$ 
if $a_i=1$ for some $i$.

$(c)_n$: $\ $ $m_n(a_1 \otimes \dots \otimes a_n)=0$ if $a_i=1$ for some $i$.

\noindent
The strictly unital condition is $(c)_n$.  The unital morphism
condition is $(a)_n$.

We first prove $(a)_3$. For $a\in H$, 
\[ 
f_2(1 \otimes a)=-Q\lambda_2(1 \otimes a)=-Q(a)=0,
\]
since $Q|_{H}=0$.  Similarly, $f_2(a \otimes 1)=0$.  This proves $(a)_3$. 
Now suppose for some $n \geq 3$ that $(a)_i$ holds for all 
$3 \leq i \leq n$.  By definition,
\[ 
\lambda_n=\sum_{s=1}^{n-1} (-1)^{s+1} \lambda_2( f_s\otimes f_{n-s}). 
\]
If $a_1=1$, $(a)_n$ implies that 
\[ 
\lambda_n(a_1 \otimes \dots \otimes a_n)=f_{n-1}(a_2 \otimes \dots
\otimes a_n)\in L. 
\]
Similarly, if $a_n=1$, we have $\lambda_n(a_1 \otimes \dots \otimes a_n)\in L$.
If $a_i=1$ for $1<i<n$, then $\lambda_n(a_1 \otimes \dots \otimes a_n)=0$. 
Therefore $(a)_i$ for $i \leq n$ implies $(b)_n$. Since $p(L)=0$,
$(c)_n$ follows from $(b)_n$. Since $Q(L)=0$, $(a)_{n+1}$ follows from
$(b)_n$.  Induction completes the proof.
\end{proof}

\begin{lemma}
\label{yylemA.5}
Let $(A,\partial)$ be a DGA and let $e\in A^0$ be an idempotent such that 
$\partial(e)=0$. Let $D=eAe$ and $C=(1-e)A+A(1-e)$.
\begin{enumerate}
\item If $HC=0$, then 
we can choose Merkulov models so that $H_{Mer}A$ is strictly isomorphic to 
$H_{Mer}D$. As a consequence
$A$ and $D$ are quasi-isomorphic as $A_\infty$-algebras.  
\item If moreover $HA$ is Adams connected, then $H^0_A$ and $H^0_D$ 
in part (a) can be chosen to contain the unit element.
\end{enumerate}  
\end{lemma}

\begin{proof} 
First of all, $D$ is a sub-DGA of $A$ with identity $e$. 
Since $A=D\oplus C$ as chain complexes,
the group of coboundaries $B^{n}$ decomposes as $B^{n} = B^{n}_{D} \oplus
B^{n}_{C}$, where $B^{n}_{D} = B^{n} \cap D$ and $B^{n}_{C} = B^{n}
\cap C$.  Since $HC=0$, we can choose $H$ and $L$ so that they
decompose similarly (with $H_{C}=0$), giving the following direct sum
decompositions: 
\begin{gather*} 
A^n=D^n\oplus C^n= (B^n_D\oplus H^n_D\oplus L^n_D)\oplus 
(B^n_C\oplus L^n_C), \\
A^n=B^n\oplus H^n\oplus L^n=(B^n_D\oplus B^n_C)\oplus H^n_D\oplus
(L^n_D\oplus L^n_C).
\end{gather*}
It follows from the construction before Theorem \ref{yythmA.2} that 
$H_{Mer}A=H_{Mer}D$. We choose $H^0_D$ to contain $e$. By Lemma 
\ref{yylemA.4}, $e$ is the strict unit of $H_{Mer}D$; hence $e$ is 
the strict unit of $H_{Mer}A$, but note that the unit $1$ of $A$ may not be 
in $HA$.  

Now suppose $HA$ is Adams connected with unit $u$. Let $H^0=ku\oplus 
H^0_{\geq 1}$ (or $H^0=ku\oplus H^0_{\leq -1}$ if negatively connected 
graded). Replace $H^0$ by $k1\oplus H^0_{\geq 1}$ and keep the 
other subspaces $B^n$, $H^n$, and $L^n$ the same. Let $\overline{H_{Mer}A}$ 
denote the new Merkulov model with the new choice of $H^0$. Then 
by Lemma \ref{yylemA.4}, $1$ is the strict unit of $\overline{H_{Mer}A}$. 
By construction, we have $(\overline{H_{Mer}A})_{\geq 1}=
(H_{Mer}A)_{\geq 1}$ as $A_\infty$-algebras without unit. By the 
unital condition, we see that $\overline{H_{Mer}A}$ is strictly isomorphic to 
$H_{Mer}A$. 
\end{proof}

\subsection{The bar construction and Ext}
\label{sectapp2}

The bar/cobar construction is one of the basic tools in homological 
algebra. Everything in this subsection is well-known -- see
\cite{FHT1}, for example -- but we 
need the details for the proof in the next subsection.

Let $A$ be a connected graded algebra and let $k$ be the trivial 
$A$-module.  Of course, the $i$-th Ext-group $\Ext^i_A(k_A,k_A)$ 
can be computed by the $i$-th cohomology of the complex 
$\Hom_A(P_A,k_A)$ where $P_A$ is any projective (or free) resolution 
of $k_A$.  Since $P_A$ is projective, $\Hom_A(P_A,k_A)$ is 
quasi-isomorphic to $\Hom_A(P_A,P_A)=\End_A(P_A)$; hence 
$\Ext^i_A(k_A,k_A)\cong H^i(\End_A(P_A))$.  Since $\End_A(P_A)$ is a 
DGA, the graded vector space $\Ext^*_A(k_A,k_A):=
\bigoplus_{i\in {\mathbb Z}} \Ext^i_A(k_A,k_A)$ has a natural algebra 
structure, and it also has an
$A_\infty$-structure by Kadeishvili's result Theorem \ref{yythmA.1}. 
By \cite[Chap.2]{Ad}, the Ext-algebra of a graded algebra $A$ 
can also be computed by using the bar construction on $A$, which 
will be explained below.

First we review the shift functor. Let $(M, \partial)$ be a complex
with differential $\partial$ of degree 1, and
let $n$ be an integer. The $n$th shift of $M$, denoted by $S^n(M)$, is
defined by
\[
S^n(M)^{i}=M^{i+n}
\]
and the differential of $S^n(M)$ is 
\[
\partial_{S^n(M)}(m)=(-1)^n \partial(m)
\]
for all $m\in M$. If $f: M\to N$ is a homomorphism of degree $p$, then 
$S^n(f): S^n(M)\to S^n(N)$ is defined by the formula
\[
S^n(f)(m)=(-1)^{pn} f(m)
\]
for all $m\in S^n(M)$. The functor $S^n$ is an automorphism of 
the category of complexes. 

The following definition is essentially standard, although sign
conventions may vary; we use the conventions from
\cite[Sect.19]{FHT1}.  Let $A$ be an augmented DGA with augmentation
$\epsilon: A\to k$, viewing $k$ as a trivial DGA. Let $I$ be the
kernel of $\epsilon$ and $SI$ be the shift of $I$. The \emph{tensor
coalgebra} on $SI$ is
\[
T(SI)=k\oplus SI\oplus (SI)^{\otimes 2}\oplus (SI)^{\otimes 3} 
\oplus \cdots,
\]
where an element $Sa_1\otimes Sa_2 \otimes \cdots \otimes Sa_n$ in 
$(SI)^{\otimes n}$ is written as
\[
[a_1| a_2| \cdots |a_n]
\]
for $a_i\in I$, together with the comultiplication
\[
\Delta([a_1 |\cdots | a_n])=
\sum_{i=0}^{n} [a_1|\cdots |a_i]\otimes [a_{i+1}|\cdots |a_{n}].
\]
The degree of $[a_1|\cdots|a_n]$ is $\sum_{i=1}^n (\deg a_i -1)$.

\begin{definition}
\label{yydefA.6} 
Let $(A,\partial_A)$ be an augmented DGA and let $I$ denote the augmentation 
ideal $\ker (A\to k)$. The {\it bar construction} on $A$ is the coaugmented 
differential graded coalgebra (DGC, for short) $BA$ defined as follows:

$\bullet$ As a coaugmented graded coalgebra, $BA$ is the tensor coalgebra 
$T(SI)$.

$\bullet$ The differential in $BA$ is the sum $d=d_0+d_1$ of the coderivations 
given by 
\[
d_0([a_1| \cdots | a_m])=-\sum_{i=1}^{m} (-1)^{n_i} [a_1 |\cdots| 
\partial_A(a_i)|\cdots | a_m]
\] 
and
\begin{gather*}
d_1([a_1])=0
\\
d_1([a_1 |\cdots |a_m])=\sum_{i=2}^{m} (-1)^{n_i} [a_1|\cdots |a_{i-1}a_i|
\cdots |a_m]
\end{gather*}
where $n_i=\sum_{j<i} (-1+\deg a_j)=\sum_{j<i} \deg [a_j]$. 
\end{definition}

The cobar construction $\Omega C$ on a coaugmented DGC $C$ is defined dually 
\cite[Sect.19]{FHT1}. We omit the definition since it is used only 
in two places, one of which is between Lemma \ref{yylemA.11}
and Lemma \ref{yylemA.12}, and the other is in Lemma \ref{yylemA.13}.

In the rest of this subsection we assume that 
$A$ is an augmented associative algebra.
In this case $SI$ is concentrated in degree $-1$; hence the degree of 
$[a_1|\cdots |a_{m}]$ is $-m$.  This means that
the bar construction $BA$ is graded by the {\it negative} of tensor 
length. The degree of the differential $d$ is $1$. We may think 
of $BA$ as a complex with $(-i)$th term equal to $I^{\otimes i}$;
the differential $d$ mapping
$I^{\otimes i}$ to $I^{\otimes i-1}$. If $A$ has an Adams grading, 
denoted $\adeg$, then $BA$ has a bigrading that is defined by
\[
\deg\; [a_1|\cdots |a_{m}]=(-m, \sum_i \adeg a_i).
\]
The second component is the Adams degree of $[a_1|\cdots |a_{m}]$.  

The bar construction on the left $A$-module $A$, denoted by $B(A,A)$, 
is constructed as follows. As a complex $B(A,A)=BA\otimes A$ with 
$(-i)$th term equal to $I^{\otimes i}\otimes A$. We use 
\[
[a_1|\cdots|a_m]x
\]
to denote an element in $I^{\otimes i}\otimes A$ where $x\in A$ and $a_i\in I$. 
The degree of $[a_1|\cdots|a_m] x$ is $-m$. The differential on $B(A,A)$ is 
defined by
\[
d(x)=0 \ \ \text{($m=0$ case)},
\]
and
\[
d([a_1 |\cdots | a_m]x)=\sum_{i=2}^{m} (-1)^{i-1} [a_1|\cdots |a_{i-1}a_i|
\cdots |a_m]x+(-1)^m [a_1|\cdots |a_{m-1}]a_m x. 
\]
Then $B(A,A)$ is a complex of free right $A$-modules.
One basic property is that the augmentations of $BA$ and $A$ make it
into a free resolution of $k_A$,
\begin{equation}
\label{E3.1.1}B(A,A)\to k_A\to 0
\end{equation}
(see \cite[19.2]{FHT1} and \cite[Chap.2]{Ad}). 

\begin{remark}
\label{yyremA.7}
In the next subsection we use the tensor $\otimes$ notation instead of the 
bar $|$ notation, which seems more natural when we concentrate on 
each term of the bar construction. 
\end{remark}

We now assume that $A$ is connected graded and finite-dimensional 
in each degree. The grading of $A$ is the Adams 
grading. Then $B(A,A)$ is bigraded with Adams grading on the second 
component, and the differential of $B(A,A)$ preserves the Adams grading.
Let $B^{\#}A$ be the graded $k$-linear dual of the coalgebra $B A$. 
Since $BA$ is locally finite, $B^{\#}A$ is a 
locally finite bigraded algebra.  With respect to the Adams grading, $B^{\#}A$ 
is {\it negatively} connected graded. The DGA $\End_A(B(A,A)_A)$ is bigraded 
too, but not Adams connected. Since $B(A,A)$ is a left differential graded 
comodule over $BA$, it has a left differential graded module structure over 
$B^{\#}A$, which is compatible with the right $A$-module structure. By an 
idea  similar to \cite[Ex. 4, p.\ 272]{FHT1} (also see \cite{LPWZ1}) 
one can show  that the natural map $B^{\#}A\to \End_A(B(A,A)_A)$ is 
a quasi-isomorphism of DGAs. 

Define the Koszul dual of a connected graded ring $A$ to be the DGA
$\End_A(P_A)$, where $P_A$ is any free resolution of $k_A$.  By the following lemma,
this definition makes sense up to quasi-isomorphism in the category of 
$A_\infty$-algebras.

\begin{lemma}
\label{yylemA.8} 
Let $A$ be a connected graded algebra and let $P_A$ and $Q_A$
be two free resolutions of $k_A$.
\begin{enumerate}
\item $\End_A(P_A)$ is quasi-isomorphic to $\End_A(Q_A)$ as 
$A_\infty$-algebras.
\item $\End_A(P_A)$ is quasi-isomorphic to $B^{\#}A$ as $A_\infty$-algebras.
\end{enumerate}
\end{lemma}

\begin{proof}
(a) We may assume that $Q_A$ is a minimal free resolution of $k_A$. Then 
$P_A=Q_A\oplus I_A$ where $I_A$ is another complex of free modules
such that $HI_A=0$ \cite[10.1.3 and 10.3.4]{AFH}. 
In this case $D:=\End_A(Q_A)$ is a sub-DGA of 
$E:=\End_A(P_A)$ such that $D=eEe$ where $e$ is the projection onto
$Q_A$. Let $C=(1-e)E+E(1-e)$. Then 
\[ 
C=\Hom_A(I_A,Q_A)+\Hom_A(Q_A,I_A)+\Hom_A(I_A,I_A),
\]
and $HC=0$. By Lemma \ref{yylemA.5}, $D$ and $E$ are quasi-isomorphic.

(b) Since $B(A,A)$ is a free resolution of $k_A$, then part (a) says
that $\End_A(P_A)$ is 
quasi-isomorphic to $\End_A(B(A,A)_A)$. The assertion follows
from the fact that $\End_A(B(A,A))$ is quasi-isomorphic to $B^{\#}A$
\cite[Ex.\ 4, p.\ 272]{FHT1}.
\end{proof}

So we may think of the bigraded DGA $B^{\#}A$ as the Koszul dual of
$A$.  This viewpoint of Koszul duality is also taken by Keller in
\cite{Ke1}.  By results in \cite{LPWZ1}, we can define the Koszul dual
of any connected graded (or augmented) $A_\infty$-algebra, and the
double Koszul dual is quasi-isomorphic to the original
$A_\infty$-algebra.

The classical Ext-algebra $\Ext^*_A(k_A,k_A)$ is the cohomology ring
of $\End_A(P_A)$, where $P_A$ is any free resolution of $k_A$. The
above lemma demonstrates the familiar fact that this is independent of
the choice of $P_A$.  Since $E:=\End_A(P_A)$ is a DGA, by Proposition
\ref{yypropA.3}, $\Ext^*_A(k_A,k_A)=HE$ has a natural
$A_\infty$-structure, which is called an {\it $A_\infty$-Ext-algebra}
of $A$. By abuse of notation we use $\Ext^*_A(k_A,k_A)$ to denote an
$A_\infty$-Ext-algebra.

\subsection{$A_\infty$-structure on Ext-algebras} 
\label{sectapp3}

In this subsection we consider the multiplications on an 
$A_\infty$-Ext-algebra of a connected graded algebra, and finally
give proofs of Theorem \ref{yythm2.2} and Proposition~\ref{yyprop2.3}.
Consider a connected graded algebra
\[
A=k\oplus A_1\oplus A_2\oplus \cdots,
\]
which is viewed as an $A_\infty$-algebra concentrated 
in degree 0, with the grading on $A$ being the Adams grading. Let $V\subset
A$ be a minimal graded vector space which generates $A$. Then 
$V\cong \fm/\fm^2$ where $\fm:=A_{\geq 1}$ is the unique maximal graded 
ideal of $A$. Let $R\subset T\langle V \rangle$ be the minimal graded 
vector space which generates the relations of $A$. Then $A\cong 
T\langle V\rangle/(R)$ where $(R)$ is the ideal generated by $R$, and
the start of a minimal graded free resolution of the trivial right $A$-module $k_A$ 
is
\begin{equation}
\label{A.8.1}
\cdots \to R\otimes A\to V\otimes A \to A\to k\to 0.
\end{equation}

\begin{lemma}
\label{yylemA.9} 
Let $A$ be a connected graded algebra. 
Then there are natural isomorphisms of graded vector spaces 
\[
\Ext^1_{A}(k_A, k_A)\cong V^{\#}
\quad \text{and}\quad \Ext^2_{A}(k_A,k_A)\cong R^{\#}.
\]
\end{lemma}

\begin{proof} 
This follows from the minimal free resolution \eqref{A.8.1}.
\end{proof}

In the rest of the section, we assume that $A$ is generated by $V=A_1$
and that $A_1$ is finite-dimensional; hence $A=T\langle A_1\rangle /(R)$.
Let $E$ be the $A_\infty$-Ext-algebra $\Ext^*_{A}(k_A,k_A)$.  We would like to 
describe the $A_\infty$-structure on $E$ by using Merkulov's
construction. 

We first fix some notation. Since $A$ is generated by $A_1$, 
the multiplication map
$\mu_n: A_1\otimes A_{n-1}\to A_n$ is surjective for $n\geq 2$. 
Let $\xi_n: A_{n}
\to A_1\otimes A_{n-1}$ be a $k$-linear map such that the composition 
\begin{equation}
\label{A.9.1}
A_n \xrightarrow{\xi_n}
A_1\otimes A_{n-1} \xrightarrow{-\mu_n} A_n
\end{equation}
is the identity map of $A_n$. Let $\theta_n$ be the composition
\[
A_n \xrightarrow{\xi_n} A_1\otimes A_{n-1}
\xrightarrow{id_{A_1}\otimes \xi_{n-1}}
A_1\otimes A_1\otimes A_{n-2} \longrightarrow \cdots \longrightarrow
A_1^{\otimes n}.
\]
Define $\xi_1=\theta_1=id_{A_1}$.  
Inductively we see that $\theta_n=(id_{A_1}\otimes \theta_{n-1})\circ \xi_n$.

Let $R=\boplus_{n\geq 2} R_{n} \subset T\langle A_1 \rangle$ be 
a minimal graded vector space of the relations of $A$.  If $A$ is a 
quadratic algebra, then $R=R_2$, but in general $R$ can be any graded
subspace of $T\langle A_1\rangle_{\geq 2}$.

We claim that we may choose $R$ such that 
\begin{equation}
\label{A.9.2}
R_n\subset A_1\otimes \theta_{n-1}(A_{n-1})\subset A_1^{\otimes n}.
\end{equation}
Inductively,
we assume that $R_{i}\subset A_1\otimes \theta_{i-1}(A_{i-1})$ for all $i<n$. 
Let $(R)_{n-1}$ be the degree $n-1$ part 
of the ideal $(R)$. Then it is generated 
by the relations of degree less than $n$, and we have a decomposition
\[
A_1^{\otimes n-1}=\theta_{n-1}(A_{n-1})\oplus \ker \mu=
\theta_{n-1}(A_{n-1})\oplus (R)_{n-1},
\]
where $\mu: A_1^{\otimes n-1}\to A_{n-1}$ is the multiplication.  
Hence we have 
\[ 
A_1^{\otimes n}=\bigl( A_1\otimes \theta_{n-1}(A_{n-1})\bigr) 
\oplus \bigl(A_1\otimes (R)_{n-1}\bigr). 
\]
Any relation $r\in R_n$ is a sum of $r_1\in A_1\otimes \theta_{n-1}(A_{n-1})$
and $r_2\in A_1\otimes (R)_{n-1}$. Modulo the relations of degree less than 
$n$, we may assume $r_2=0$. Hence the claim is proved. 

Assuming now that $R_n\subset A_1\otimes \theta_{n-1}(A_{n-1})\subset
\fm\otimes \fm$ for all $n$, then the minimal resolution \eqref{A.8.1}
becomes a direct summand of the bar resolution \eqref{E3.1.1}
\[ 
\cdots \to \fm^{\otimes 2}\otimes A\to \fm\otimes A \to A\to k\to 0. 
\]
The proof of Theorem \ref{yythm2.2} follows from Merkulov's construction 
and several lemmas below.

\begin{remark}
\label{yyremA.10}
After we identify $T(S\fm)$ with a complex, all linear maps
between subspaces/dual spaces of $(S\fm)^{\otimes i}$ now have degree 0, so no
extra sign will be introduced when commuting two linear maps. 
\end{remark}

Since $A$ is concentrated in degree 0, the grading on the differential 
graded coalgebra $T(S\fm)$ is by the negative of the wordlength, namely, 
$(T(S\fm))^{-i}=\fm^{\otimes i}$. The differential $d=(d^i)$ of the bar 
construction $T(S\fm)$ is induced by the multiplication 
$\fm \otimes \fm \to \fm$ in $A$. For example, 
\[
d^{-1}([a_1])=0 \quad \text{and}\quad d^{-2}([a_1|a_2])=(-1)^{-1}[a_1 a_2]
\]
for all $a_1, a_2\in \fm$. 
There is a natural decomposition of $\fm$ with respect to the Adams
grading,
\[
\fm=A_1\oplus A_2\oplus A_3\oplus \cdots,
\]
which gives rise to a decomposition of $\fm\otimes \fm$ with respect to the 
Adams grading:
\[
\fm\otimes \fm=(A_1\otimes A_1)\oplus (A_1\otimes A_2\oplus
A_2\otimes A_1) \oplus \cdots.
\]
Using the map
\[ 
R_n\to A_1\otimes \theta_{n-1}(A_{n-1})\to A_1\otimes A_{n-1}, 
\]
we can view $R_n$ as a subspace of $A_1\otimes A_{n-1}$. 

\begin{lemma}
\label{yylemA.11} Let $W^{-2}_{n}=\boplus_{i+j=n} A_i\otimes A_j$, where $n$
is the Adams grading. Then there are decompositions of vector spaces
\begin{gather*}
W^{-2}_{n}=\im (d^{-3}_{n})\oplus R_n\oplus \xi_n(A_n), \\
\ker (d^{-2}_{n})=\im (d^{-3}_{n}) \oplus R_n,
\end{gather*}
where $R_n$ and $\xi_n(A_n)$ are subspaces of $A_1\otimes A_{n-1}$.
\end{lemma}

\begin{proof} It is clear that the injection 
$A_n \xrightarrow{\xi_n} A_1\otimes A_{n-1}
\longrightarrow W^{-2}_{n}$ defines a projection from 
$W^{-2}_{n}$ to $A_n$. Since $d^{-2}_{n}: 
W^{-2}_{n}\to A_n$ is a surjection,
we have a decomposition
\[
W^{-2}_{n}=\ker (d^{-2}_{n})\oplus \xi_n(A_n).
\]
Since $R^{\#}\cong \Ext^2_{A}(k_A,k_A)=H^2((T(S\fm))^{\#})$
by Lemma \ref{yylemA.9},
there is a decomposition $\ker (d^{-2}_{n})=\im (d^{-3}_{n}) 
\oplus R_n $ where
$R_n\subset A_1\otimes A_{n-1}$ by the choice of $R$. 
Hence the assertion follows.
\end{proof}

This lemma says that $\im (d_{n}^{-3})$ is in the ideal generated 
by the relations
of degree less than $n$. 

Since $A$ is Adams locally finite, $(\fm^{\otimes n})^{\#}\cong
(\fm^{\#})^{\otimes n}$ for all $n$. Let $\Omega A^{\#}$ be the cobar 
construction on the DGC $A^{\#}$. Via the isomorphisms
\[
B^{\#}A=(T(S\fm))^{\#}\cong T((S\fm)^{\#})\cong T(S^{-1}\fm^{\#})
=\Omega A^{\#},
\]
we identify $B^{\#}A=(T(S\fm))^{\#}$ with $\Omega A^{\#}=T(S^{-1}\fm^{\#})$.
The differential $\partial$ on $B^{\#}A$ is defined by
\[
\partial(f)=-(-1)^{\deg f} f\circ d
\]
for all $f\in T(S^{-1}\fm^{\#})$ .


We now study the first two nonzero differential maps of $\Omega A^{\#}$,
\[
\partial^1: \fm^{\#}\to (\fm^{\#})^{\otimes 2}\qquad\text{and}\qquad 
\partial^2: (\fm^{\#})^{\otimes 2}\to  (\fm^{\#})^{\otimes 3}.
\] 
For all $s$ and $n$, let
\[
T^s=(\fm^{\#})^{\otimes s}
\]
and
\[
T^s_{-n}=\boplus_{i_1+\cdots +i_s=n} A^{\#}_{i_1}\otimes \cdots \otimes
A^{\#}_{i_s}.
\]
Fix the Adams degree $-n$, and let 
\begin{gather*}
\partial^{1}_{-n}:A_n^{\#}\to \boplus_{i+j=n} A_i^{\#}\otimes A_j^{\#}, \\
\partial^{2}_{-n}:\boplus_{i+j=n} A_i^{\#}\otimes A_j^{\#}\to 
\boplus_{i_1+i_2+i_3=n}A_{i_1}^{\#}\otimes A_{i_2}^{\#}\otimes A_{i_3}^{\#}.
\end{gather*}
The decomposition \eqref{A.1.1} for $T^1_{-n}$ is
\begin{gather*}
T^1_{-1}=H^1_{-1}, \quad\text{ i.e., } \quad B^1_{-1}=L^1_{-1}=0, \\
T^1_{-n}=L^1_{-n}\quad \text{ i.e., }\quad B^1_{-n}=H^1_{-n}=0,
\end{gather*}
for
all $n\geq 2$. The decomposition \eqref{A.1.1} for $T^2_{-n}$ is given
below. 

\begin{lemma} 
\label{yylemA.12} Fix $n \geq 2$.  With notation as above, we have the
following.
\begin{enumerate}
\item Define the duals of subspaces by using the decompositions given in 
Lemma \ref{yylemA.11}. Then $\im  \partial^1_{-n}=(\xi_n(A_n))^{\#}$ and 
$\ker \partial^2_{-n} =(\xi_n(A_n))^{\#}\oplus R_n^{\#}.$
\item 
The decomposition \eqref{A.1.1} for $T^{2}_{-n}$ can be chosen to be
\[
T^2_{-n}=B^2_{-n}\oplus H^2_{-n} \oplus L^2_{-n}=
(\xi_n(A_n))^{\#}\oplus R_n^{\#}\oplus (\im d^{-3}_{n})^{\#}.
\] 
The projections onto $R_n^{\#}$ and $(\xi_n(A_n))^{\#}$ kill $\boplus_{i+j=n, j>1}
A_i^{\#}\otimes A_j^{\#}$.
\item Let $Q$ be the homotopy defined in Merkulov's construction for
DGA $T(\fm^{\#})$. Then we may choose $Q^{2}_{-n}$ to be equal to
$-(\xi_n)^{\#}$ when 
restricted to $T^2_{-n}$.
\end{enumerate}
\end{lemma} 

\begin{proof} (a) This follows from Lemma \ref{yylemA.11} and a linear 
algebra argument. 

(b) This follows from Lemma \ref{yylemA.11}, part (a), and the fact
that $R_n$ and
$\xi_n(A_n)$ are subspaces of $A_1\otimes A_{n-1}$. 

(c) Let $\xi_n$ also denote the map $A_n\to A_1\otimes A_{n-1}\to W^{-2}_n$.
Since $d^{-2}_n=-\mu_n$, the composite $d^{-2}_n \circ \xi_n: 
A_n\to W^{-2}_n\to A_n$ is the identity map (see \eqref{A.9.1}). 
Since $ (d^{-2}_n \circ 
\xi_n)^{\#}=(\xi_n)^{\#} \circ (d^{-2}_n)^{\#}$, 
\[
(\xi_n)^{\#}\circ (d^{-2}_n)^{\#}:\quad A_n^{\#}\to T^2_{-n}\to A_n^{\#}
\]
is the identity map of $A_n^{\#}$, which is a subspace $T^1$. 
Since $(d^{-2}_n)^{\#}=-\partial^1_{-n}$, we may 
choose the homotopy $Q$ to be $-(\xi_n)^{\#}$ when restricted to $T^2_{-n}$.
\end{proof}

When $n=1$, we formally set $Q\lambda_1=-id_T$. 
Let $\lambda_i$ be defined as in \eqref{A.1.2} in Merkulov's construction.
In particular, $\lambda_2$ is the multiplication of $T(S^{-1}\fm^{\#})$.

Recall that $T(S^{-1}\fm^{\#})$ is a free (or tensor) DGA generated by 
$S^{-1}\fm^{\#}$. Then for all $a_1\otimes\cdots\otimes a_n\in 
(S^{-1}\fm^{\#})^{\otimes n}$ and $b_1\otimes \cdots \otimes b_m
\in (S^{-1}\fm^{\#})^{\otimes m}$ we have
\begin{equation}
\label{A.12.1}
\lambda_2((a_1\otimes\cdots\otimes a_n)\botimes (b_1\otimes \cdots
\otimes b_m))=(a_1\otimes \cdots\otimes a_n\otimes b_1\otimes 
\cdots\otimes b_m).
\end{equation}
By the above formula, we see that $\lambda_2$ changes $\botimes$ to $\otimes$,
so it is like the identity map. 

\begin{lemma} 
\label{yylemA.13} 
Let $E^1=\Ext^1_{A}(k_A,k_A)=A_1^{\#}$. Let $n\geq 2$.
\begin{enumerate}
\item When restricted to $(E^1)^{\otimes n}$, the map $\lambda_n$  has
image in 
\[
T^2_{-n} =\boplus_{i+j=n} A_i^{\#}\otimes A_j^{\#}.
\]
\item When restricted to $(E^1)^{\otimes n}$, the map $-Q_n\lambda_n$ is 
the $k$-linear dual of the map
\[
\theta_n: A_n\to A_1\otimes A_{n-1}\to A_1^{\otimes n}.
\]
\item When restricted to $(E^1)^{\otimes n}$, the map $m_n=Pr_H \lambda_n$ is the
$k$-linear  the dual of the canonical map
\[
R_n\longrightarrow A_1\otimes A_{n-1}
{\overset{id\otimes \theta_{n-1}}{\longrightarrow}} (A_1)^{\otimes n}.
\]
\end{enumerate}
\end{lemma} 

\begin{proof} We use induction on $n$. 

(a)  By definition,
\[
\lambda_n=\lambda_2 \sum_{\substack{i+j=n, \\ i,j>0}} (-1)^{i+1}
Q\lambda_i \botimes Q \lambda_j.
\]
For $n=2$, the claim follows from \eqref{A.12.1}. Now assume $n>2$. 
By (b) each $Q\lambda_i$ has image in $(A_i)^{\#}$ for all $i<n$. Hence 
the image of $\lambda_n$ is in 
$\boplus_{i+j=n} (A_i)^{\#}\otimes (A_j)^{\#}$. 

(b) When $n=2$, $\theta_2=\xi_2$, and the claim follows from Lemma 
\ref{yylemA.12}(c). 
Now we assume $n>2$. When restricted to $(E^1)^{\otimes n}$, the image of
$Q\lambda_i\botimes Q\lambda_j$ is in $A_i^{\#}\botimes A_j^{\#}$ for all
$i+j=n$. If $i>1$, 
\[
Q\lambda_2(A_i^{\#}\botimes A_j^{\#})=Q(A_i^{\#}\otimes A_j^{\#})=
-(\xi_n)^{\#}(A_i^{\#}\otimes A_j^{\#})=0
\]
(see Lemma \ref{yylemA.12}(b,c)).
Therefore, when restricted to $(E^1)^{\otimes n}$, we have 
\[
Q\lambda_n=Q\lambda_2[(-1)^{2} (-id)\botimes Q\lambda_{n-1}]=
-Q\lambda_2(id\botimes Q\lambda_{n-1}).
\]
By induction, $Q\lambda_{n-1}= -(\theta_{n-1})^{\#}$, and by Lemma 
\ref{yylemA.12}(c) we see that $Q=-(\xi_n)^{\#}$.  Hence 
\[
Q\lambda_n=(\xi_n)^{\#} \lambda_2 (id_{A^{\#}_1}\botimes 
-(\theta_{n-1})^{\#})
=-((id_{A_1}\otimes (\theta_{n-1})))\circ \xi_n)^{\#}=-(\theta_n)^{\#}.
\]

(c) Since we assume that each $R_n$ is subspace of $A_1\otimes A_{n-1}$, 
then the dual of 
\[
R_n\to A_1\otimes A_{n-1}
\]
is $Pr_H$ restricted to $A_1^{\#}\otimes A_{n-1}^{\#}$. Hence the dual of
$R_n\to A_{1}^{\otimes n}$ is equal to $Pr_H\circ 
(id\otimes \theta_{n-1})^{\#}$.

By Lemma \ref{yylemA.12}(b), $Pr_H$ is zero when applied to 
$(A_i)^{\#}\otimes (A_j)^{\#}$ for all  $i>1$. By (a),
\[
Pr_H \lambda_n =Pr_H \lambda_2 (-id_{A^{\#}_1}\otimes Q\lambda_{n-1})
=Pr_H (id_{A^{\#}_1}\otimes (\theta_{n-1})^{\#}),
\]
which is the desired map.
\end{proof} 

\begin{proof}[Proof of Theorem \ref{yythm2.2}] First of all by 
(A.9.2) and its proof, we may assume that $R_{n}\subset A_1\otimes 
\theta_{n-1}(A_{n-1})$. So Lemmas \ref{yylemA.11}, \ref{yylemA.12} and
\ref{yylemA.13} hold. The canonical map in Lemma \ref{yylemA.13}(c)
is just the inclusion. Therefore the assertion holds.
\end{proof}

The following corollary is immediate.

\begin{corollary} 
\label{yycorA.14} Let $A$ and $E$ be as in Theorem \ref{yythm2.2}. 
\begin{enumerate}
\item The algebra $A$ is determined by the 
maps $m_n$ restricted to $(E^1)^{\otimes n}$ for all $n$.
\item The $A_\infty$-structure of $E$ is determined up to 
quasi-isomorphism by the 
maps $m_n$ restricted to $(E^1)^{\otimes n}$ for all $n$.
\end{enumerate}
\end{corollary}

\begin{proof}
(a) By Theorem \ref{yythm2.2}, the map 
\[
R\to \boplus_{n} (A_1)^{\otimes n}
\] can be 
recovered from $m_n$ restricted to $(E^1)^{\otimes n}$. Hence the structure
of $A$ is determined. 

(b) After $A$ is recovered, the $A_\infty$-structure of $E$ is determined 
by $A$. Therefore the structure of $E$ is determined by the 
restriction of $m_n$ on $(E^1)^{\otimes n}$, up to quasi-isomorphism.
\end{proof}

\begin{proof}[Proof of Proposition \ref{yyprop2.3}]
By the construction given above, it is clear that if the grading group
for the Adams grading 
is $\Z\oplus G$ for some abelian group $G$, then all of the maps 
including $m_n$ preserve the $G$-grading.  The assertion follows.
\end{proof}

{\it Acknowledgments} $\quad$ D.-M. Lu is supported by 
the Pao Yu-Kong and Pao Zhao-Long Scholarship. Q.-S. Wu is 
supported by the NSFC (project 10171016, key project 
10331030), in part by the grant of STCSM:03JC14013 in China, 
and by the Ky/Yu-Fen Fan Fund of the American Mathematical 
Society. J.J. Zhang is supported by NSF grant DMS-0245420 
(USA) and Leverhulme Research Interchange Grant F/00158/X (UK). 
A part of research was done when J.J. Zhang was visiting the 
Institute of Mathematics, Fudan University, China. J.J. Zhang thanks 
Fudan University for the warm hospitality and the financial support.



\begin{thebibliography}{LPWZ2}


\bibitem[Ad]{Ad}
J.F. Adams, 
\emph{On the non-existence of elements of {H}opf invariant one}, Ann.
of Math. (2) \textbf{72} (1960), no.1, 20--104.

\bibitem[An]{An}
D. Anick, 
\emph{On the homology of associative algebras}, Trans. Amer. Math. Soc. 
\textbf{296} (1986), no. 2, 641--659.

\bibitem[ASZ]{ASZ}
M. Artin, L.W. Small, and J.J. Zhang, 
\emph{Generic flatness for strongly Noetherian algebras,} 
J. Algebra \textbf{221} (1999), no. 2, 579--610.

\bibitem[ASc]{ASc} 
M. Artin and W. Schelter, \emph{Graded algebras of global dimension $3$}, 
Adv. in Math. \textbf{66} (1987), no. 2, 171--216. 

\bibitem[ASt1]{ASt1}
M. Artin and J.T. Stafford, 
\emph{Noncommutative graded domains with quadratic 
growth}, Invent. Math. \textbf{122} (1995), no. 2, 231--276.

\bibitem[ASt2]{ASt2}
M. Artin and J.T. Stafford, \emph{Semiprime graded algebras of dimension two},
J. Algebra \textbf{227} (2000), no. 1, 68--123.
 
\bibitem[ATV1]{ATV1} 
M. Artin, J. Tate, and M. Van den Bergh, 
\emph{Some algebras associated to automorphisms of elliptic curves}, 
The Grothendieck Festschrift,
Vol. I, 33--85, Progr. Math., 86, Birkh{\" a}user Boston, Boston, MA, 1990.

\bibitem[ATV2]{ATV2} 
M. Artin, J. Tate, and M. Van den Bergh, 
\emph{Modules over regular algebras of dimension $3$}, Invent. Math. 
\textbf{106} (1991), no. 2, 335--388. 

\bibitem[AZ]{AZ}
M. Artin and J.J. Zhang, \emph{Noncommutative projective schemes}, 
Adv. Math. \textbf{109} (1994), no. 2, 228--287.

\bibitem[AFH]{AFH} 
L.L. Avramov, H.-B. Foxby, and S. Halperin,
Differential graded homological algebra, manuscript. 





\bibitem[Be]{Be}
G. M. Bergman, 
\emph{The diamond lemma for ring theory}, 
Adv. in Math. \textbf{29} (1978), no. 2, 178--218.

\bibitem[Ek]{Ek}
E.K. Ekstr{\" o}m,
\emph{The Auslander condition on graded and filtered Noetherian rings,} 
S{\' e}minaire d'Alg{\` e}bre Paul Dubreil et Marie-Paul Malliavin, 
39{\` e}me Ann{\' e}e (Paris, 1987/1988), 220--245, 
Lecture Notes in Math., 1404, Springer, Berlin, 1989. 

\bibitem[FHT]{FHT1}
Y. F{\'e}lix, S. Halperin, and J.-C. Thomas, \emph{Rational homotopy theory},
Springer-Verlag, New York, 2001.







\bibitem[Ka]{Ka}
T.V. Kadeishvili,  \emph{On the theory of homology of fiber spaces}, 
(Russian) International Topology Conference (Moscow State Univ., Moscow, 
1979). Uspekhi Mat. Nauk \textbf{35} (1980), no. 3(213), 183--188.
Translated in Russ. Math. Surv. \textbf{35} (1980), no.~3, 231--238.



\bibitem[Ke1]{Ke1}
B. Keller, \emph{Deriving {D}{G} categories}, Ann. Sci. \'Ecole Norm.
Sup. (4) \textbf{27} (1994), no.~1, 63--102.

\bibitem[Ke2]{Ke2}
\bysame, \emph{Bimodule complexes via strong homotopy actions}, Algebr.
Represent. Theory \textbf{3} (2000), no.~4, 357--376, Special issue dedicated
to Klaus Roggenkamp on the occasion of his 60th birthday.

\bibitem[Ke3]{Ke3} 
\bysame, \emph{Introduction to ${A}$-infinity algebras and modules}, Homology
Homotopy Appl. \textbf{3} (2001), no.~1, 1--35 (electronic).

\bibitem[Ke4]{Ke4}
\bysame, 
\emph{$A$-infinity algebras in representation theory,}  
Contribution to the Proceedings of ICRA IX, Beijing 2000. 

\bibitem[Ke5]{Ke5} 
\bysame,
\emph{$A_\infty$-algebras and triangulated categories}, 
in preparation.

\bibitem[Ko]{Ko} M. Kontsevich, \emph{Homological algebra of mirror 
symmetry}, Proceedings of the International Congress of
Mathematicians, Vol. 1, 2 (Z\"urich, 1994), 120--139, Birkh{\"{a}}user, 
Basel, 1995.


\bibitem[LBS]{LBS}
L. Le Bruyn and S.P. Smith, \emph{Homogenized ${\mathfrak sl}_2$},
Proc. Amer. Math. Soc. \textbf{118} (1993), no. 3, 725--730. 

\bibitem[LBSV]{LBSV}
L. Le Bruyn, S.P. Smith, and M. Van den Bergh,
\emph{Central extensions of three-dimensional Artin-Schelter regular
algebras}, Math. Z. \textbf{222} (1996), no. 2, 171--212.


\bibitem[Lev]{Lev} 
T. Levasseur, \emph{Some properties of noncommutative
regular rings}, Glasgow Math. J. \textbf{34} (1992), 277-300.

\bibitem[LS]{LS}
T. Levasseur and S.P. Smith,\emph{Modules over the $4$-dimensional Sklyanin
algebra}, Bull. Soc. Math. France \textbf{121} (1993), no. 1, 35--90.

\bibitem[LPWZ1]{LPWZ2}
D.-M. Lu, J.H. Palmieri, Q.-S. Wu, and J.J. Zhang, 
\emph{$A_\infty$-algebras for ring theorists}, Algebra Colloquium,
\textbf{11} (2004), no. 1, 91-128.

\bibitem[LPWZ2]{LPWZ1}
\bysame, \emph{Koszul equivalences 
in $A_\infty$-algebras}, in preparation.

\bibitem[LWZ]{LWZ}
D.-M. Lu, Q.-S. Wu, and J.J. Zhang, \emph{Regular algebras of type 
(1,2,2,2,1) and their Ext-algebras},  in preparation. 


\bibitem[MR]{MR} 
J. C. McConnell and J. C. Robson, 
\emph{Noncommutative Noetherian Rings,} Wiley, Chichester, 1987.

\bibitem[Me]{Me1}
S.A. Merkulov, \emph{Strong homotopy algebras of a K{\" a}hler manifold}, 
Internat. Math. Res. Notices 1999, no. \textbf{3}, 153--164. 




\bibitem[SW]{SW} M. Shi and J. Wang, \emph{$A_\infty$-algebras and AS 
regular algebras of dimension 3}, preprint, (2003).

\bibitem[Sk1]{Sk1}
E.K. Sklyanin, \emph{Some algebraic structures connected with the 
Yang-Baxter equation} (Russian), Funktsional. Anal. i Prilozhen
\textbf{16} (1982), no. 4, 27--34.

\bibitem[Sk2]{Sk2}
E.K. Sklyanin, \emph{Some algebraic structures connected with the
Yang-Baxter equation. Representations of a quantum algebra}, (Russian) 
Funktsional. Anal. i Prilozhen \textbf{17} (1983), no. 4, 34--48.

\bibitem[SV1]{SV1}
B. Shelton and M. Vancliff,
\emph{Some Quantum ${\mathbb P}^3$s with One Point,}
Comm. Alg. \textbf{27} No. 3 (1999), 1429-1443. 

\bibitem[SV2]{SV2}
B. Shelton and M. Vancliff,
\emph{Embedding a Quantum Rank Three Quadric in a Quantum ${\mathbb P}^3$,}
Comm. Alg. \textbf{27} No. 6 (1999), 2877-2904.

\bibitem[SSW]{SSW}
L.W. Small, J.T. Stafford, and R.B. Warfield, 
\emph{Affine algebras of Gel\'fand-Kirillov dimension one are PI}, 
Math. Proc. Cambridge Philos. Soc. \textbf{97} (1985), no. 3, 407--414.


\bibitem[Sm]{Sm1}
S. P. Smith, 
\emph{Some finite-dimensional algebras related to elliptic curves,}
Representation theory of algebras and related topics (Mexico City, 1994), 
315--348, CMS Conf. Proc., \textbf{19}, AMS, Providence, RI, 1996. 

\bibitem[SS]{SS} 
S.P. Smith and J.T. Stafford, 
\emph{Regularity of the four-dimensional Sklyanin algebra,}
Compositio Math. \textbf{83} (1992), no. 3, 259--289. 

\bibitem[Staf]{Staf} 
J.T. Stafford, 
\emph{Regularity of algebras related to the Sklyanin algebra}, 
Trans. Amer. Math. Soc. \textbf{341} (1994), no. 2, 895--916.


\bibitem[Sta]{Sta} J. D. Stasheff, 
\emph{Homotopy associativity of $H$-spaces. I, II}, 
Trans. Amer. Math. Soc. \textbf{108} (1963), 275-292; ibid.
\textbf{108} (1963) 293--312. 

\bibitem[Ste1]{Ste1} D.R. Stephenson, 
\emph{Artin-Schelter regular algebras of global dimension three}, 
J. Algebra \textbf{183} (1996), no. 1, 55--73.

\bibitem[Ste2]{Ste2} D.R. Stephenson, 
\emph{Algebras associated to elliptic curves}, 
Trans. Amer. Math. Soc. \textbf{349} (1997), no. 6, 2317--2340.

\bibitem[Ste3]{Ste3} D.R. Stephenson, 
\emph{Quantum planes of weight $(1,1,n)$}, 
J. Algebra \textbf{225} (2000), no. 1, 70--92. 

\bibitem[SteZ]{SZ}
D.R. Stephenson and J.J. Zhang,
\emph{Growth of graded Noetherian rings}, 
Proc. Amer. Math. Soc. \textbf{125} (1997), no. 6, 1593--1605. 

\bibitem[VV1]{VV1}
K. Van Rompay and M. Vancliff,
\emph{Embedding a Quantum Nonsingular Quadric in a Quantum ${\mathbb P}^3$,}
J. Algebra \textbf{195} No. 1 (1997), 93-129.

\bibitem[VV2]{VV2}
K. Van Rompay and M. Vancliff,
\emph{Four-dimensional Regular Algebras with Point Scheme a Nonsingular Quadric
in ${\mathbb P}^3$}, Comm. Alg. \textbf{28} No. 5 (2000), 2211-2242.

\bibitem[VVW]{VVW}
K. Van Rompay, M. Vancliff, and L. Willaert,
\emph{Some Quantum ${\mathbb P}^3$s with Finitely Many Points},
Comm. Alg. \textbf{26} No. 4 (1998), 1193-1208.

\bibitem[Va1]{Va1}
M. Vancliff, 
\emph{Quadratic Algebras Associated with the Union of a Quadric 
and a Line in ${\mathbb P}^3$}, J. Algebra \textbf{165} No. 1 (1994),
63-90.

\bibitem[Va2]{Va2}
M. Vancliff, 
\emph{The Defining Relations of Quantum n x n Matrices}, 
J. London Math. Soc. \textbf{52} No. 2 (1995), 255-262.

\bibitem[Zh]{Zh1}
J.J. Zhang, 
\emph{Twisted graded algebras and equivalences of graded 
categories}, Proc. London Math. Soc. (3) \textbf{72} (1996), 
no. 2, 281--311. 
 

\end{thebibliography}
\end{document}